\documentclass[12pt]{article}

\usepackage{amsmath, amsthm, amscd, amsfonts, amssymb}

 \newtheorem{theorem}{Theorem}[section]

 \newtheorem{lemma}[theorem]{Lemma}
 \newtheorem{definition}[theorem]{Definition}
 
 \newtheorem{remark}[theorem]{Remark}

  \usepackage[usenames,dvipsnames]{color}
  
\def\divv{{\rm div}}

\def\rrd{{\mathbb{R}^d}}
\def\call{{\mathcal{L}}}

\def\calf{{\mathcal{F}}}

\def\calo{{\mathcal{O}}}

\def\cald{{\mathcal{D}}}
\def\calx{{\mathcal{X}}}

\def\call{{\mathcal{L}}}

\def\vsp{\vspace*{1,5mm}\\ }

\def\mk{\medskip }
\def\sk{\smallskip }
\def\n{\noindent }
\def\dd{\displaystyle}
\def\barr{\begin{array}}
\def\earr{\end{array}}
\def\FP{Fokker--Planck}

\def\1{^{-1}}

\def\rr{{\mathbb{R}}}
\def\nn{{\mathbb{N}}}
\def\9{{\infty}}
\def\lbb{{\lambda}}
\def\wt{\widetilde}
\def\ov{\overline}
\def\vf{{\varphi}}

\def\vp{{\varepsilon}}

\def\ff{\forall }
\def\({\left(}
\def\){\right)}
\def\<{\left<}
\def\>{\right>}

 \title{Nonlocal, nonlinear Fokker--Planck equations and nonlinear martingale problems} 
 \author{Viorel Barbu\thanks{Al.I. Cuza University and Octav Mayer Institute of Mathematics of  Romanian Academy, Ia\c si, Romania.  Email: vbarbu41@gmail.com}\and Jos\'e Lu\'\i s da Silva\thanks{CIMA -- Faculty of Exact Sciences and Engineering, University of Madeira, 9020-105 Funchal, Portugal. Email: joses@staff.uma.pt}\and Michael R\"ockner\thanks{Fakult\"at f\"ur Mathematik, Universit\"at Bielefeld,  D-33501 Bielefeld, Germany. Email: roeckner@math.uni-bielefeld.de}\ \thanks{Academy of Mathematics and System Sciences, CAS, Beijing, China.}}
 \date{}

\begin{document}
	
\maketitle

\begin{abstract}
This work is concerned with the existence  of mild solutions and the uniqueness of distributional solutions to non\-linear Fokker--Planck equations with nonlocal operators $\Psi(-\Delta)$, where $\Psi$ is a Bernstein function. As applications, the existence and uniqueness of solutions to the corresponding nonlinear martingale problems are proved. Furthermore, it is shown that these solutions form a nonlinear Markov process in the sense of McKean such that their one-dimensional time marginal law densities are the solutions to the nonlocal nonlinear \FP\ equation. Hence, McKean's program envisioned in his PNAS paper from 1966 is realized for these nonlocal PDEs.\sk\\
{\bf Keywords:}  Fokker--Planck  equation,  Bernstein function, distributional solution, mild solution, nonlinear martingale problem.
\end{abstract}

 \section{Introduction}\label{s1}
We are  concerned here with the following nonlinear, nonlocal \FP\ equation of Nemytskii-type
\begin{equation}\label{e1.1}
	\barr{l}
	u_t+\Psi(-\Delta)(\beta(u))+{\rm div}(Db(u)u)=0,\ \mbox{ on  } (0,\9)\times\rr^d,\\
	u(0,x)=u_0(x),\ x\in\rr^d,
	\earr\end{equation}\newpage\n where the functions $\beta:\rr\to\rr$, $D:\rrd\to\rrd$, $d\ge2$,   $b:\rr\to\rr$ are  to be made precise below, and $\Psi$ is a Bernstein function. 

The operator $\Psi(-\Delta)$ is defined as follows. Let $S':=S'(\rr^d)$ be the dual of the Schwartz test function space $S:=S(\rr^d)$. Define 
\begin{equation}\label{e1.1p}
	D_{\Psi}:=\{u\in S':\calf(u)\in L^1_{\rm loc}(\rr^d),\Psi(|\xi|^2)\calf(u)\in S'\}(\supset L^1(\rr^d))\end{equation}
and $\Psi(-\Delta):D_\Psi\to S'$ by
\begin{equation}\label{e1.2p}
	\calf(\Psi(-\Delta)u)(\xi):=
	\Psi(|\xi|^2)\calf(u)(\xi),\ \xi\in\rr^d,\ u\in D_{\Psi}, \end{equation}
where $\calf$ stands for the Fourier transform on $\rr^d$, i.e., 
$$\calf(u)(\xi)=(2\pi)^{-d/2}
\int_{\rr^d}e^{ix\cdot\xi}u(x)dx,\ \xi\in\rr^d,\ u\in L^1(\rr^d).$$($\calf$ extends from $S'$ to itself.)

Furthermore, $\Psi:[0,\9)\to[0,\9)$ is a Bernstein function, i.e. an infi\-nitely differentiable completely monotone function, that is,
$$(-1)^k\Psi^{(k)}(r)\ge0,\ \ff r\ge0,\ k=1,2,...\,.$$A Bernstein function $\Psi$ admits the unique representation (see \cite{12''}, p.~21)
\begin{equation}\label{e12}
	\Psi(r)=a_1+a_2r+\int^\9_0(1-e^{-rt})\mu(dt),\ \ff r\ge0,
\end{equation}
where $a_1,a_2\ge0$ and $\mu$ is a positive measure on $(0,\9)$ such that
\begin{equation}\label{e1.4'}
	m:=\int^\9_0(1\wedge t)\mu(dt)<\9,\vspace*{-2mm}
\end{equation}which implies 
	$\Psi(r)\le m(1+r),\ r\ge0.$

Given a Bernstein function $\Psi$, there is a unique convolution semigroup of sub-probability (if $a_1=0$, probability) measures $(\eta^\Psi_t)_{t\ge0}$ on $(0,\9)$ such that
\begin{equation}\label{e1.4a}
	\call(\eta^\Psi_t)(\lbb)=e^{-t\Psi(\lbb)},\ \ff\lbb\ge 0,\end{equation}where $\call(\eta^\Psi_t)$ is the Laplace transform of $\eta^\Psi_t$ (see \cite{12''}, p.~49).

A standard example is $\Psi(r)\equiv r^s$, $s\in(0,1)$, which corresponds to \mbox{$a_1,a_2=0$} and
\begin{equation}\label{e1.7a}
	\mu(dt)=\frac s{\Gamma(1-s)}\ t^{-s-1}dt.\end{equation}
The hypotheses below will be in effect in Sections \ref{s2} and \ref{s4.1}. 
\begin{itemize}
	\item[\rm(i)] $\beta\in C^1(\rr)\cap {\rm Lip}(\rr),\ \beta(0)=0,\ \beta'(r)>0,\ \ff\,r\in\rr.$\vspace*{-2mm} 
	\item[\rm(ii)] $D\in (C^1\cap C_b)(\rrd;\rrd),\ {\rm div}\,D\in L^2_{\rm loc}(\rrd).$\vspace*{-2mm}  
	\item[\rm(iii)] $b\in C_b(\rr)\cap C^1(\rr).$\vspace*{-2mm} 
	\item[\rm(iv)] $({\rm div}\,D)^-\in L^\9,\ b(r)\ge0,\ \ff r\in\rr$,  where $({\rm div}\,D)^-=-\inf((\divv\,D,0))$.\vspace*{-2mm} 
	\item[(v)] $\Psi:[0,\9)\to[0,\9)$ is a Bernstein function of the form \eqref{e12} 
	with $a_1,a_2=0$ which satisfies for some $s\in\(\frac12,1\)$ and $C\in(0,\9)$,\vspace*{-2mm}
	\begin{equation}\label{e1.4aa}
		\Psi(r)\ge Cr^s,\ \ff r\ge0.\end{equation}
\end{itemize}
Here, we shall study the existence of a {\it mild solution} to equation \eqref{e1.1} (see Definition \ref{d1.1} below) and also the uniqueness of distributional solutions. As~regards the existence, we shall follow the semigroup methods developed in \mbox{\cite{a2''}--\cite{2}} (see also \cite{6b}) and in the special case $\Psi(r)\equiv r^s,$ $s\in\(\frac12,1\)$ in \cite{a9} (see also \cite{6b}). Namely, we shall represent  \eqref{e1.1} as an abstract differential equation in $L^1(\rrd)$ of the form\vspace*{-2mm}	
\begin{equation}\label{e1.4}
	\dd\frac{du}{dt}+A(u)=0,\  t\ge0,\ \ 	u(0)=u_0,\end{equation}where $A:D(A)\subset L^1(\rrd)\to L^1(\rrd)$ is an $m$-accretive realization in $L^1(\rrd)$ of the operator\vspace*{-1mm} 
\begin{equation}\label{e1.5}
\barr{rcl}
	A_0(u)&\!\!\!=\!\!\!&\Psi(-\Delta)(\beta(u))+{\rm div}(Db(u)u),\  u\in D(A_0),\sk\\
	D(A_0)&\!\!\!=\!\!\!&\left\{u\in L^1(\rrd);\Psi(-\Delta)(\beta(u))+{\rm div}(Db(u)u)\in L^1(\rrd)\right\}.\earr\end{equation}	

\begin{definition}\label{d1.1} \rm A function $u\in C([0,\9);L^1(\rrd))$ is said to be a {\it mild solution} to \eqref{e1.1} with initial condition $u_0\in L^1$, if for each $0<T<\9$,   
	\begin{eqnarray}
		&u(t)=\dd\lim_{h\to0}u_h(t)\mbox{ in }L^1(\rrd)\mbox{ uniformly in }t\in[0,T], \label{e1.6}\end{eqnarray}where $u_h$ is for each $h>0$ the solution to the equation\vspace*{-1mm} 
	\begin{equation}\label{e16}
		\barr{ll}
		\dd\frac1h\ (u_h(t)-u_h(t-h))+A_0 u_h(t)=0,&\ff t\ge0,\sk\\ 
		u_h(t)=u_0,&\ff t\le0.\earr\end{equation}
\end{definition}	
We note that, if $u$ is a mild solution to \eqref{e1.1}, then it is also a {\it Schwartz distributional solution}, that is,  
\begin{eqnarray}
	 	&&\int^\9_0\dd\int_\rrd (u(t,x)
	\vf_t(t,x)-\Psi(-\Delta)(\vf(t,x))\beta(u(t,x))\nonumber\\[-1mm]
	&&\qquad\quad\qquad\quad+\,b(u(t,x))u(t,x)D(x)\cdot\nabla\vf(t,x))dtdx\label{e1.10}\vspace*{-2mm} \\ 
&&\qquad\quad\qquad\quad+\int_\rrd \vf(0,x)u_0(dx)=0,\ \ff\,\vf\in C^\9_0([0,\9)\times\rrd),\nonumber
\end{eqnarray}
where, in general,  $u_0$ is a measure of finite variation on $\rrd$. 
The existence  for equation \eqref{e1.1} is given by Theorem \ref{t2.3}.  
In Section \ref{s3}, the uniqueness of dis\-tri\-bu\-tio\-nal solutions to \eqref{e1.1}, \eqref{e1.10} respectively, in the class\break \mbox{$(L^1\cap L^\9)((0,T)\times\rrd)$}  $\cap L^\9(0,T;L^2)$ will be proved under Hypotheses (j)--(jjj) on $D,b$ and $\beta$ stated in Section~\ref{s3}. 
We would like to stress at this point that as in \cite{4} (see, also, the pio\-nee\-ring papers \cite{5} and \cite{a15'}), where such a result was proved for local nonlinear \FP\ equations (i.e. of type \eqref{e1.1} with $-\Delta$ replacing $\Psi(-\Delta)$), we prove uniqueness in the large class of distributional solutions without any a priori-restrictions such as e.g. weak differentiability of the solutions.   
Therefore, our uniqueness results are considerably stronger than uniqueness results within the much smaller classes of mild solutions or in the local case of entropy solutions (see, e.g., \cite{8a}). At this point, we would also like to draw attention to the paper \cite{9prim}, where uniqueness of distributional solutions for nonlocal porous media equations was also proved, but in the case where $D=0$ and by a technique quite different from ours.  
In Section \ref{s4} we apply the above  analytic results on \eqref{e1.1} to prove exis\-tence and uniqueness (see Theorems \ref{t4.1} and \ref{t4.3}, resp.) for the nonlinear martingale problem corresponding to \eqref{e1.1} in the spirit of \cite{a11'}.  
Furthermore, we prove that the solutions to such a nonlinear martingale problem form a nonlinear Markov process in the sense of McKean \cite{17a} (see Theorem \ref{t4.4} below). The existence proof (Theorem \ref{t4.1}) is based on Theorem \ref{t2.3} and \cite{7aa}, in which a nonlocal analogue of the superposition principle (see \cite{b17',24'}) is proved. The uniqueness proof (Theorem \ref{t4.3}) is based on Theorems \ref{t3.1} and our result on linearized uniqueness of \eqref{e1.1} (see Theorem \ref{t3.2}). 
As a consequence of this and Corollary 3.8 in \cite{18a}, we finally obtain Theorem \ref{t4.4}, which in particular realizes McKean's vision from \cite{17a}, i.e. to identify the solutions of non\-linear PDEs as the one-dimensional time marginal law densities of a non\-li\-near Markov process, for nonlocal PDEs as in \eqref{e1.1}.
In the special case $\Psi(r)\equiv r^s$, $s\in\(\frac12,1\)$, equation \eqref{e1.1} was  studied  \cite{a9} (see, also, \cite{7,7a,11'} for a direct approach to existence theory in the case $\Psi(r)=r^s$ and $D\equiv0$). Though the strategy in the present paper to prove the existence and uniqueness of solutions to \eqref{e1.1} follows in great lines the one  developed in \cite{a9}. In fact, we deliberately keep the structure of the proofs close to the respective ones in \cite{a9} so that the partly quite substantial differences between \cite{a9} and this paper become clearly visible. For example, the crucial inequality \eqref{e2.15a} in this paper is much more general than the correspon\-ding inequality (also \eqref{e2.15a}) in \cite{a9}.   Furthermore, we have to work in  special Sobolev spaces defined in terms of $\psi$ (see the beginning of Section \ref{s2}). In addition, the proof of the uniqueness Theorem \ref{t3.1} is different from that of the corresponding result in \cite{a9} (i.e., Theorem 3.1 in there) in one of its essential parts. Moreover, the probabilistic applications in Section \ref{s4} are more difficult to prove in the case of general Bernstein functions compared to the case $\psi(r)=r^s$, $r\ge0$, $s\in\(\frac12,1\)$. We also like to point out that the proofs in \cite{a9} were written in a very concise way, leaving too many quite nontrivial parts to the reader. Here, we fill in these parts to a reasonable extent, so that this paper may be also considered as completing \cite{a9}, as far as such parts are concerned. 

Finally, we would like to mention that the present paper represents a significant improvement compared to its first arXiv version.

\mk\noindent{\bf Notation.} $L^p(\rrd)=L^p,\ p\in[1,\9]$ is the standard space of Lebesgue  \mbox{$p$-integrable} functions on $\rr^d$. We denote the correspon\-ding local space by $L^p_{\rm loc}$  and  the norm of $L^p$ by $|\cdot|_p$. The inner product in $L^2$ is denoted by $(\cdot,\cdot)_2$.  Let $H^\sigma(\rrd)=H^\sigma$, $0<\sigma<\9$, denote the standard Bessel space on $\rrd$ and  $H^{-\sigma}$ its dual space. Let $C_b(\rr)$ denote the space of continuous and bounded functions on $\rr$ and   $C^1(\rr)$ the space of continuously dif\-fe\-ren\-tiable functions on~$\rr$, and likewise $C^1(\rrd,\rrd)$ the space of continuously differentiable vector fields from $\rrd$ to $\rrd$. For any $T>0$ and a Banach space $\calx$, $C([0,T];\calx)$ denotes the space of $\calx$-valued continuous functions on $[0,T]$ and by $L^p(0,T;\calx)$ the space of $\calx$-valued $L^p$-Bochner integrable functions on $(0,T).$   $C^k_c(\calo)$, $\calo\subset\rrd$, denotes the space of  $k$-differentiable functions with compact support in $\calo$ and   $\cald'(\calo)$  the space of Schwartz distributions on $\calo$.   $C^\9_0([0,\9)\times\rrd)$   denotes the space of dif\-fe\-ren\-tiable functions on $[0,\9)\times\rrd$ with compact in $[0,\9)\times\rrd$.   \mbox{$S'(\rrd)=S'$}   denotes the space of tempered distributions on $\rrd$ and   $\calf(y)$ the Fourier transform of $y\in S'(\rrd)$.  

\section{Existence of a mild solution}\label{s2}
\setcounter{equation}{0}

We first note that by \eqref{e1.4'}  
all functions
$$\rrd\ni\xi\to(\Psi(\vp+|\xi|^2))^\alpha,\ \vp>0,\ \alpha\in\rr,$$are multipliers on $S$, hence on $S'$. 

Hence, we may define the maps $\Psi(\vp I-\Delta):S'\to S'$ by
$$\Psi(\vp I-\Delta)u:=\calf\1(\Psi(\vp+|\xi|^2)\calf u),\ u\in S',$$which are clearly linear homeomorphisms, and the following Hilbert spaces:
\newpage
$$\barr{c}
H^\Psi:=H^\Psi(\rr^d):=\left\{u\in S':\sqrt{\Psi(1+|\xi|^2)}\,\calf(u)\in L^2\right\},\vsp
\dot H^\Psi:=\dot H^\Psi(\rr^d):=\left\{u\in S':\calf(u)\in L^1_{\rm loc}\mbox{ and }\sqrt{\Psi(|\xi|^2)}\,\calf(u)\in L^2\right\},\earr$$with respective norms
\begin{eqnarray*}
|u|^2_{H^\Psi}&:=&\int_{\rr^d}\Psi(1+|\xi|^2)|\calf(u)(\xi)|^2d\xi\vsp
|u|^2_{\dot H^\Psi}&:=&\int_{\rr^d}\Psi(|\xi|^2)|\calf(u)(\xi)|^2d\xi.\end{eqnarray*}
We denote the corresponding inner products by $\<\cdot,\cdot\>_{H^\Psi}$ and $\<\cdot,\cdot\>_{\dot H^\Psi}$, respectively.  
By \eqref{e1.4'}  and \eqref{e1.4aa} we have the continuous embeddings
\begin{equation}\label{2.0}
	H^1\subset H^\Psi\subset H^s\subset\dot H^s,\ \ H^\Psi\subset \dot H^\Psi\subset\dot H^s,\end{equation}where $\dot H^s$ denotes the usual homogeneous Sobolev space of order $s$. We note that $\dot H^s$ is only complete, if $s<\frac d2$, which holds in our case since $s<1$, $d\ge2$. Clearly, $H^\Psi$ is complete, but since by \eqref{e1.4aa}
$$\frac1{\Psi(|\xi|^2)}\le\frac1C\ \frac1{|\xi|^{2s}}\in L^1_{\rm loc},$$it follows that also $\dot H^\Psi$ is complete. For $\vp\ge0$ we define
$$D_{\vp,2}:=\{u\in L^2:\Psi(\vp I-\Delta)u\in L^2\}.$$Then, it is elementary to check that $(\Psi(\vp I-\Delta),D_{\vp,2})$ is a nonnegative definite self-adjoint operator on $L^2$ and that, for all $u\in D_{0,2}$,
$$|u|^2_{\dot H^\Psi}=\int_\rrd\left|\sqrt{\Psi(-\Delta)}\,u\right|^2d\xi.$$Furthermore, $D_{\vp,2}:=D_{1,2},\ \ff\vp>0$, and  
$$|u|^2_{H^{\Psi,\vp}}=\int_\rrd \left|\sqrt{\Psi(\vp I-\Delta)}\,u\right|^2d\xi,$$where $|\cdot|_{H^{\Psi,\vp}}$ is defined analogously to $|\cdot|_{H^\Psi}$ with $\Psi(\vp+|\xi|^2)$    replacing $\Psi(1+|\xi|^2).$ Then, we have the following

\begin{lemma}\label{l2.0}\
	\begin{itemize}
		\item[\rm(i)] $\dot H^\Psi$ and $H^{\Psi,\vp}$, $\vp>0$, as well as $H^\alpha$, $\alpha\in(0,1]$,  are invariant under composition with Lipschitz continuous functions $\vf:\rr\to\rr$ with \mbox{$\vf(0)=0$.} 
		\item[\rm(ii)] For all $\alpha\in\rr$ and $\vp>0$, we have for the inverse $\Psi(\vp I-\Delta)\1:S'\to S'$ of $\Psi(\vp I-\Delta):S'\to S'$ that
		$$\Psi(\vp I-\Delta)\1(H^\alpha)\subset H^{\alpha+2s}$$and that the operator 
		$\Psi(\vp I-\Delta)\1:H^\alpha\to H^{\alpha+2s}$ is continuous.
	\end{itemize}
\end{lemma}

\begin{proof} (i): This is an immediate consequence of \cite[Section 1.5, in particular, Theorem 1.5.3 and Example 1.5.2]{16'}, since, as it is proved there, $\(\dot H^\Psi,\<\cdot,\cdot\>_{\dot H^\Psi}\)$, $\(\dot H^{\Psi,\vp},\<\cdot,\cdot\>_{\dot H^{\Psi,\vp}}\)$, $\vp>0$, and $H^\alpha$, $\alpha\in(0,1]$ are transient Dirichlet spaces. 

(ii): Using the definition of $H^\alpha,\ \alpha\in\rr$, in terms of Fourier transforms,   
the proof is elementary  by~\eqref{e1.4aa}.\end{proof} 

We shall now prove  the following key lemma which is similar to Lemma 2.1 in \cite{a9}.  (See, also, \cite{a2''}--\cite{2}, \cite{6b}.)

\begin{lemma}\label{l2.1}  Define
	\begin{equation}\label{e2.1prim}
	\lbb_0:=(|(\divv\,D)^-|_\9|b|_\9)^{-1},\end{equation}where we set $\frac10:=\9$. 
	Then, under Hy\-po\-theses {\rm(i)--(v)} there is a family of operators  $\{J_\lbb:L^1\to L^1;\lbb>0)\}$, which for all $\lbb\in(0,\lbb_0)$  satisfies
	\begin{eqnarray}
		&	(I+\lbb A_0)(J_\lbb(f))=f,\ \ff\,f\in L^1,\label{e2.1}\\ 
		&	|J_\lbb(f_1)-J_\lbb(f_2)|_1\le|f_1-f_2|_1,\ \ff f_1,f_2\in L^1,\label{e2.2}\\ 
		&J_{\lbb_2}(f)=J_{\lbb_1}\(\dd\frac{\lbb_1}{\lbb_2}\,f+\(1-\frac{\lbb_1}{\lbb_2}\)J_{\lbb_2}(f)\),\ \ff f\in L^1,\ \lbb_1,\lbb_2\in(0,\lbb_0),\qquad\label{e2.3}\\ 
		&\dd\int_\rrd J_\lbb(f)dx=\int_\rrd f\,dx,\ \ff f\in L^1,\label{e2.4}\\ 
		& J_\lbb(f)\ge0,\ \mbox{ a.e. on }\rrd,\mbox{ if }f\ge0,\mbox{ a.e. on }\rrd,\label{e2.5}\\ 
		& |J_\lbb(f)|_\9\le
		\(1+\dd\frac\lbb{\lbb_0}\)
		|f|_\9,\ \ff\,f\in L^1\cap L^\9,	\label{e2.5a}\\ 
		&\beta(J_\lbb(f))\in H^{\Psi}\cap L^1\cap L^\9,\ \ff f\in L^1\cap L^\9. \label{e2.6a}
	\end{eqnarray}
 Furthermore, $J_\lbb(f)$ is the unique solution in $D(A_0)\cap L^1\cap L^\9$ to \eqref{e2.1}, if $f\in L^1\cap L^\9$.
\end{lemma}

\begin{proof} 
	We shall  first prove the existence  solutions for the equation
\begin{equation}\label{e2.6} 
	y+\lbb A_0(y)=f\mbox{ in }S',\end{equation}where $f\in L^1$.  We consider the approximating equation
\begin{equation}\label{e2.7} 
	y+\lbb\Psi (\vp  I-\Delta)(\beta_\vp(y))+\lbb\,{\rm div}(D_\vp b_\vp(y)y)=f\mbox{ in }S',\end{equation}where $\vp\in(0,1]$, $\beta_\vp(r):=\beta(r)+\vp r$ and
	$$D_\vp:=\eta_\vp D,\ \eta_\vp(x):=\eta(\vp|x|),\,x\in\rrd,\mbox{ with }\eta\in C^1_b([0,\9)),\,{\bf1}_{[0,1]}\le\eta\le{\bf1}_{[0,2]}.$$
We have
\begin{equation}
	\label{e2.8a}
	\barr{c}
	|D_\vp|\in L^2\cap L^\9,\ |D_\vp|\le|D|,\ \dd\lim_{\vp\to0}D_\vp(x)=D(x),\mbox{ a.e. }x\in \rrd,\vsp 
	{\rm div}\,D_\vp\in L^2,\ ({\rm div}\,D_\vp)^-\le({\rm div}\,D)^-+
	\vp|\eta'|_\9|D|,\earr\end{equation}
where we used that $\nabla\eta_\vp(x)=\vp\eta'(\vp|x|)\frac x{|x|},\ x\in\rrd$. Furthermore, the function $b_\vp$ is defined by
$$b_\vp(r)\equiv \frac{(b*\vf_\vp)(r)}{1+\vp|r|},\ \ff\,r\in\rr,$$where $*$ is the convolution product and $\vf_\vp(r)=\frac1\vp\ \vf\(\frac r\vp\)$ is a standard mollifier.  We set $b^*_\vp(r):=b_\vp(r)r,$ \mbox{$ r\in\rr,$} and note that   $b^*_\vp$ is bounded, Lipschitz.

Let us first assume  that $f\in L^2$ and rewrite \eqref{e2.7} as  
\begin{equation}\label{2.11}
	F_{\vp,\lbb}(y)=f\mbox{ in }S',\end{equation}where $F_{\vp,\lbb}:L^2\to S'$ is defined by
$$F_{\vp,\lbb}(y):=y+\lbb\Psi(\vp I-\Delta)\beta_\vp(y)+\lbb\,{\rm div}(D_\vp b^*_\vp(y)),\ \ff y\in L^2.$$
We set 
 $$G_\vp(y):=\Psi(\vp I-\Delta)(y),\ y\in S'.$$  
Now, we shall show that \eqref{2.11} has a unique solution $y_\vp\in L^2$. To this end, we rewrite it as  
$$G^{-1}_\vp (F_{\vp,\lbb}(y))=G^{-1}_\vp f  (\in H^{2s} \mbox{ by Lemma \ref{l2.0} (ii)}),$$that is,
\begin{equation}
	\label{2.12}
	G^{-1}_\vp y+\lbb\beta_\vp(y)+\lbb G^{-1}_\vp\,{\rm div}(D_\vp b^*_\vp(y))=G^{-1}_\vp f.
\end{equation}
Since $D_\vp b^*_\vp(y)\in L^2$, we have ${\rm div}(D_\vp b^*_\vp(y))\in H\1$, and so, by Lemma \ref{l2.0} (ii) and because $s>\frac12$, we have 
that $G^{-1}_\vp F_{\vp,\lbb}:L^2\to L^2$ is continuous. Now, it is easy to see that \eqref{2.12} has a unique solution, $y_\vp\in L^2$ for small enough $\lbb$, because, by \eqref{2.12} we have, for $y_1,y_2\in L^2$, 
\begin{eqnarray}
 &&\hspace*{-6mm}	(G^{-1}_\vp(F_{\vp,\lbb}(y_2)-F_{\vp,\lbb}(y_1)),y_2-y_1)_2 
 =(G^{-1}_\vp(y_2-y_1),y_2-y_1)_2\label{2.13}\\
 && 
	+\lbb(\beta_\vp(y_2)-\beta_\vp(y_1),y_2-y_1)_2 \nonumber\\
	&& 
	-\lbb{\ }_{H^{-1}}\!\!\<{\rm div}(D_\vp(b^*_\vp(y_2)-b^*_\vp(y_1))),
	G^{-1}_\vp(y_2-y_1)\>_{H^1}\nonumber\\
	&& \ge|y_2-y_1|^2_{H^{-s}}+\lbb\vp|y_2-y_1|^2_2 
 -\lbb c_1|D_\vp(b^*_\vp(y_2)-b^*_\vp(y_1))|_2|\nabla G^{-1}_\vp(y_2-y_1)|_2\nonumber\\ 
 &&  \ge|y_2-y_1|^2_{H^{-s}}+\lbb\vp|y_2-y_1|^2_2-\lbb c_\vp |D|_\9{\rm Lip}(b^*_\vp)|y_2-y_1|_2|y_2-y_1|_{H^{1-2s}},
	\nonumber\end{eqnarray} 
where $c_\vp\in(0,\9)$ is independent of $\lbb,y_1,y_2$.  Since $-s<1-2s<0$, by interpolation  we have for $\theta:=\frac{2s-1}s$ that
$$|y_2-y_1|_{H^{1-2s}}\le|y_2-y_1|^{1-\theta}_2|y_2-y_1|^\theta_{H^{-s}}\,,$$and so, by Young's inequality we find that the left hand side of \eqref{2.13} domi\-nates
$$\lbb(\vp-\lbb c_\vp)|y_2-y_1|^2_2+\frac12|y_2-y_1|^2_{H^{-s}}$$for some $c_\vp\in(0,\9)$ independent of $\lbb,y_1$ and $y_2$. Hence, for some $\lbb_\vp\in(0,\9)$, we conclude that $G^{-1}_\vp F_{\vp,\lbb}$ is strictly monotone on $L^2$ for $\lbb\in(0,\lbb_\vp)$.

By Lemma \ref{l2.0} (ii), it follows from \eqref{2.12} that $\beta_\vp(y_\vp)\in H^{2s-1}$, and so,  by Lemma \ref{l2.0} (i), because $s<1$ and $\beta^{-1}_\vp$ is Lipschitz and zero at zero, the solution $y_\vp$ also belongs to $H^{2s-1}$, hence $b^*_\vp(y_\vp)\in H^{2s-1}$. Since $s\in\(\frac12,1\)$ and $D\in C^1(\rrd;\rrd)$, by simple bootstrapping \eqref{2.12} it follows that 
$y_\vp\in H^1,$ 
and therefore,  by \eqref{2.12} we have 
$	\beta_\vp(y_\vp)\in G^{-1}_\vp (L^2)$ and so 
\begin{equation}\label{e2.13''}
G_\vp\beta_\vp(y_\vp)\in L^2.\end{equation} Furthermore, for $f\in L^2$ and $\lbb\in(0,\lbb_\vp)$, $y_\vp$ is the unique solution to \eqref{2.11}  in $L^2$.   

Assume now that $\lbb\in(0,\lbb_\vp)$ and that  $f\ge0$, a.e. on $\rrd$. 
Then, we have
\begin{equation}
	\label{e2.13}
	y_\vp\ge0,\ \mbox{ a.e. on }\rrd.
\end{equation}
To prove this, consider the function

\begin{equation}
	\label{e2.14}
	\eta_\delta(r)=\left\{\barr{rll}
	-1&\mbox{ for }&r\le-\delta,\vsp
	\dd\frac r\delta&&r\in(-\delta,0),\vsp 
	0&\mbox{ for }&r\ge0,\earr\right.
\end{equation}
$\delta>0$, and multiply  \eqref{e2.7}, where $y=y_\vp$, by $\eta_\delta(\beta_\vp(y_\vp))$ $(\in H^1)$ and integrate over $\rrd$. By \eqref{e2.13''} we get  
\begin{equation}
	\label{e2.15}
	\barr{l}
	\dd\int_\rrd y_\vp\eta_\delta(\beta_\vp(y_\vp))dx
	+\lbb\dd\int_\rrd G_\vp(\beta_\vp(y_\vp))\eta_\delta(\beta_\vp(y_\vp))dx\\
	\qquad 
	=\dd\int_\rrd f\eta_\delta(\beta_\vp(y_\vp))dx
	+\lbb\dd\int_\rrd D_\vp b^*_\vp(y_\vp)\eta'_\delta(\beta_\vp(y_\vp))\cdot\nabla \beta_\vp(y_\vp) dx.
	\earr\end{equation}On the other hand, we have the following crucial inequality
\begin{equation}
	\label{e2.15a}
	\hspace*{-4mm}\barr{r}
	\dd\int_\rrd G_\vp (u)\vf(u)dx
	=\<u,\vf(u)\>_{H^{\Psi,\vp}}\ge\frac12({\rm Lip}(\vf))\1\<\vf(u),\vf(u)\>_{H^{\Psi,\vp}}\ge0,\\ \ff u\in G^{-1}_\vp (L^2)(\subset H^\Psi),\earr\hspace*{-4mm}
\end{equation}for all non constant, nondecreasing Lipschitz functions $\vf:\rr\to\rr.$ This is a well-known inequality in the theory of Dirichlet forms (see, e.g., \cite[Examples 6.4 and 6.5]{23'}).   By \eqref{e2.13''}, this yields
\begin{equation}\label{e2.15aa}
	\dd\int_\rrd G_\vp(\beta_\vp(y_\vp))\eta_\delta(\beta_\vp(y_\vp))dx
	\ge 0.\end{equation}
Taking into account that $|\beta_\vp(y_\vp)|\ge \vp|y_\vp|$, we have 
\begin{eqnarray}
	\label{e2.16}
	\left|\dd\int_\rrd D_\vp b_\vp(y_\vp)y_\vp
	\eta'_\delta(\beta_\vp(y_\vp))\nabla \beta_\vp(y_\vp) dx\right|
 \le\dd\frac1\delta\,|b|_\9\int_{\wt E^\delta_\vp}|y_\vp|\,|\nabla \beta_\vp(y_\vp)|\,|D_\vp|dx\ \ \\[-3mm]
\le\dd\frac 1\vp\,|b|_\9\|D_\vp\|_{L^2}
	\(\dd\int_{\wt E^\delta_\vp}|\nabla \beta_\vp(y_\vp)|^2dx\)^{\frac12}\to0\mbox{ as }\delta\to0.\nonumber
\end{eqnarray}
Here  $\wt E^\delta_\vp=\{-\delta<\beta_\vp(y_\vp)\le0\}$ and we used that $\nabla\beta_\vp(y_\vp)=0$, a.e. on $\{\beta_\vp(y_\vp)=0\}$.  

Since ${\rm sign}\,\beta_\vp(r)\equiv{\rm sign}\,r$, by \eqref{e2.15}--\eqref{e2.16}  we get, for $\delta\to0$, that $y^-_\vp=0,$ a.e. on $\rrd$ and so \eqref{e2.13} holds. 

If $\lbb\in(0,\lbb_\vp)$, $f\in L^2$, and $y_\vp=y_\vp(\lbb,f)$ is the solution to \eqref{e2.7} in $L^2$, we have for $f_1,f_2\in L^1\cap L^2$
\begin{eqnarray}
	\label{e2.18}
	y_\vp(\lbb,f_1)-y_\vp(\lbb,f_2)
	+\lbb G_\vp(\beta_\vp(y_\vp(\lbb,f_1))-\beta_\vp(y_\vp(\lbb,f_2)))\\
	+\lbb\ {\rm div}\ D_\vp(b^*_\vp(y_\vp(\lbb,f_1))-b^*_\vp(y_\vp(\lbb,f_2)))=f_1-f_2.\nonumber \end{eqnarray}

Now, for $\delta>0$ consider the function
\begin{equation*}
	\calx_\delta(r)=\left\{\barr{rll}
	1&\mbox{ for }&r\ge\delta,\vsp
	\dd\frac r\delta&\mbox{ for }&|r|<\delta,\vsp 
	-1&\mbox{ for }&r<-\delta.\earr\right.
\end{equation*}
Multiplying \eqref{e2.18} by $\calx_\delta(\beta_\vp(y_\vp(\lbb,f_1))-\beta_\vp(y_\vp(\lbb,f_2)))$, we get
$$\barr{l}
\dd\int_\rrd(y_\vp(\lbb,f_1)-y_\vp(\lbb,f_2))\calx_\delta(\beta_\vp(y_\vp(\lbb,f_1))-\beta_\vp(y_\vp(\lbb,f_2)))dx\\ 
\quad\le\!\lbb\,\dd\frac1\delta\!\int_{E_\vp^\delta}\!|b^*_\vp(y_\vp(\lbb,f_1)){-}b^*_\vp(y_\vp(\lbb,f_2))|\,|D_\vp|\\
\qquad\qquad\qquad
|\nabla(\beta_\vp(y_\vp(\lbb,f_1))
{-}\beta_\vp(y_\vp(\lbb,f_2)))|dx +|f_1-f_2|_1,\earr$$because,  by virtue of \eqref{e2.15a},  
$$\dd\int_{\rrd}G_\vp(\beta_\vp(y_\vp,f_1)-\beta_\vp(y_\vp,f_2))\calx_\delta(\beta_\vp(y_\vp,f_1)-\beta_\vp(y_\vp,f_2))dx\ge0.$$ 

Set $E^\delta_\vp=\{|\beta_\vp(y_\vp(\lbb,f_1))-\beta_\vp(y_\vp(\lbb,f_2))|\le\delta\}.$  

Since $|\beta_\vp(r_1)-\beta_\vp(r_2)|\ge\vp|r_1-r_2|,\ D_\vp\in L^2(\rrd;\rrd)$, $b^*_\vp\in {\rm Lip}(\rr),$ $y_\vp(\lbb,f_i)\in H^1,\ i=1,2$, and 
$\nabla(\beta_\vp(y_\vp(\lbb,f_1))-\beta_\vp(y_\vp(\lbb,f_2)))=0$, a.e. on $\{\beta_\vp(y_\vp(\lbb,f_1))-\beta_\vp(y_\vp(\lbb,f_2))=0\}$, we~get that
\begin{equation}\label{e2.22prim}
\barr{r}
\dd\lim_{\delta\to0}\frac1\delta \int_{E^\delta_\vp} |b^*_\vp(y_\vp(\lbb,f_1)){-}b^*_\vp(y_\vp(\lbb,f_2))|\,|D_\vp|\cdot\qquad\vsp
\qquad\cdot|\nabla(\beta_\vp(y_\vp(\lbb,f_1))-\beta_\vp(y_\vp(\lbb,f_2)))|dx=0.\earr\end{equation}
This yields
\begin{eqnarray}
&|y_\vp(\lbb,f_1)-y_\vp(\lbb,f_2)|_1\le|f_1-f_2|_1,\ \ff\lbb\in(0,\lbb_\vp),\label{e2.19}\\[1mm]
&
|y_\vp(\lbb,f)|_1\le|f|_1,\ \ff\,f\in L^1\cap L^2,\ \lbb\in(0,\lbb_\vp).	\label{e2.23a}
\end{eqnarray}
Now, let us remove the restriction on $\lbb\in(0,\lbb_\vp)$. To this purpose  define the operator \mbox{$A_\vp:D_0(A_\vp)\!\to\! L^1$} by
\begin{eqnarray*}
A_\vp(y)&:=&G_\vp(\beta_\vp(y))+{\rm div}(D_\vp b^*_\vp(y)),\vsp 
D_0(A_\vp)&:=&\{y\in L^1: G_\vp(\beta_\vp(y))+{\rm div}(D_\vp b^*_\vp(y))\in L^1\},\end{eqnarray*}and
$$J^\vp_\lbb(f):=y_\vp(\lbb,f),\ f\in L^1\cap L^2,\ \lbb\in(0,\lbb_\vp).$$
Then, by \eqref{2.11}, \eqref{e2.13''},  \eqref{e2.23a} and since $y_\vp(\lbb,f)\in H^1,$
\begin{equation}\label{e2.22'}
	\barr{c}
	J^\vp_\lbb(f)\in D_0(A_\vp)\cap H^1, \mbox{ and } 
 \beta_\vp(J^\vp_\lbb(f))\in G^{-1}_\vp(L^2),\vsp \ff f\in L^1\cap L^2,\ \lbb\in(0,\lbb_\vp).\earr
\end{equation} 
For $\lbb_1,\lbb_2\in(0,\lbb_\vp),\ f\in L^1\cap L^2$, we have
\begin{equation}\label{e2.25prim}
J^\vp_{\lbb_2}f+\lbb_2A_\vp(J^\vp_{\lbb_2}\,f)=f\end{equation}
and thus
$$J^\vp_{\lbb_2}f+\lbb_1A_\vp(J^\vp_{\lbb_2}\,f)=
\(1-\frac{\lbb_1}{\lbb_2}\)
	J^\vp_{\lbb_2}f+\frac{\lbb_1}{\lbb_2}\,f.$$
	Since, as seen above, $J^\vp_\lbb\,f,\lbb\in(0,\lbb_\vp)$, is the unique solution in $L^2$ of \eqref{2.11}, we obtain that
\begin{equation}\label{e2.25secund}
J^\vp_{\lbb_2}f=J^\vp_{\lbb_1}
\(\(1-\frac{\lbb_1}{\lbb_2}\)
J^\vp_{\lbb_2}f+\frac{\lbb_1}{\lbb_2}\,f\),\ 
f\in L^1\cap L^2;\ \lbb_1,\lbb_2\in(0,\lbb_\vp).
\end{equation}	 
Furthermore, by \eqref{e2.19}, $J^\vp_\lbb$ extends by continuity to an operator \mbox{$J^\vp_\lbb:L^1\to L^1$}  and \eqref{e2.25secund} extends to all $f\in L^1$.  We note that the operator $(A_\vp,D_0(A_\vp))$ is closed as an operator on $L^1$. Hence, by \eqref{e2.19},  
$J^\vp_\lbb(L^1)\subset D_0(A_\vp)$ and so,  $J^\vp_\lbb(f)$ solves \eqref{e2.7} for all $f\in L^1$.  But it might not be the unique solution in $D_0(A_\vp)$ of \eqref{e2.7}. Therefore, we define $D(A_\vp):=J^\vp_\lbb(L^1)$, $\lbb\in(0,\lbb_\vp)$, and restrict $A_\vp$ to $D(A_\vp)$. By \eqref{e2.25secund}, extended to all $f\in L^1$, it follows  that $D(A_\vp)$ is independent of $\lbb\in(0,\lbb_\vp)$  and that $(I+\lbb A_\vp):D(A_\vp)\to L^1$ is a bijection, hence $J^\vp_\lbb=(I+\lbb A_\vp)\1$. In particular,   $J^\vp_\lbb(f)$ is the unique solution in $D(A_\vp)$ of \eqref{e2.7} for all $f\in L^1,$ $\lbb\in(0,\lbb_\vp)$.

Now let $0<\lbb_1<\lbb_\vp$. Then, for $\lbb\ge\lbb_\vp$  as above, the equation
\begin{equation}
	\label{e2.23aaaa}
	y+\lbb A_\vp(y)=f  \ (\in L^1),\ y\in D(A_\vp),
\end{equation}can be rewritten as
$$y+\lbb_1 A_\vp(y)=\(1-\frac{\lbb_1}{\lbb}\)y+\frac{\lbb_1}\lbb\,f,$$
 hence 
\begin{equation}
	\label{e2.23v}
	y=J^\vp_{\lbb_1}\(\(1-\frac{\lbb_1}\lbb\)y+\frac{\lbb_1}\lbb\,f\). 
\end{equation}
Taking into account that, by \eqref{e2.19}, 
$|J^\vp_{\lbb_1}(f_1)-J^\vp_{\lbb_1}(f_2)|_1\le|f_1-f_2|_1,$ it follows that \eqref{e2.23v} has a unique solution $y_\vp\in J^\vp_{\lbb_1}(L^1)=D(A_\vp)$. Let  $J^\vp_\lbb(f):=y_\vp$, $\lbb\in[\lbb_\vp,\9),\ f\in L^1$. Then, $J^\vp_\lbb(f)$ is the unique solution in $D(A_\vp)$ of \eqref{e2.23aaaa}  for all $\lbb>0,$ $f\in L^1$, hence it is independent of $\lbb_1$ and, by the same arguments as above, it follows that
	\begin{equation}\label{e2.27prim}
	J^\vp_{\lbb_2}\,f=	J^\vp_{\lbb_1}
	\(\(1-\frac{\lbb_1}{\lbb_2}\)	J^\vp_{\lbb_2}\,f+\frac{\lbb_1}{\lbb_2}\,f\),\ f\in L^1;\ \lbb_1,\lbb_2\in(0,\9).
	\end{equation}
	By \eqref{e2.23v} we also see that \eqref{e2.19}, \eqref{e2.23a} extend to all $\lbb>0,$  $f\in L^1.$ 

Let us prove that, for $f\in L^1\cap L^2$,
\begin{equation}\label{b2.27'}
	J^\vp_\lbb(f)\in H^1\ \mbox{ and }\ \beta_\vp(J^\vp_\lbb(f))\in G^{-1}_\vp(L^2)  \mbox{ for all } \lbb>0.\end{equation}
Here is the argument. Fix $\lbb_1\in[\lbb_\vp/2,\lbb_\vp)$ and set $\lbb{:=}2\lbb_1.$ Define 
$$T(y){:=}J^\vp_{\lbb_1}\(\frac12 y+\frac12 f\),\ y\in L^1.$$
Then, as just proved, for any $f_0\in L^1\cap L^2 $ we have 	$\lim\limits_{n\to\9}T^n(f_0)=J^\vp_\lbb(f)\mbox{ in }L^1.$ It suffices to prove
\begin{equation}\label{b2.27'''}
	J^\vp_\lbb(f)\in L^2  ,\end{equation}
because then $J^\vp_\lbb(f)=J^\vp_{\lbb_1}(g)$ with $g:=\frac12\, J^\vp_{\lbb}(f)+\frac12\ f\in L^1\cap L^2,$ and so \eqref{b2.27'} follows by \eqref{e2.22'}.
 
To prove \eqref{b2.27'''}, we note that by \eqref{e2.22'} we have 
$$(I+\lbb_1 A_\vp)T^n(f_0)=\frac12\ T^{n-1}(f_0)+\frac12\, f\ \mbox{ for }n\in\nn,$$with $T^n(f_0)\in H^1$ and $\beta_\vp(T^n(f_0))\in G^{-1}_\vp(L^2).$ Hence, multiplying this equation, by $T^n(f_0)$ and integrating over $\rrd$ we find
\begin{eqnarray}\label{b2.27iv}
	&&\hspace*{-14mm}	|T^{n}f_0|^2_2+\lbb_1
	(G_\vp  (\beta_\vp(T^n(f_0))),\beta^{-1}_\vp(\beta_\vp(T^n(f_0))))_2\\     
	&&\hspace*{-10mm} =\lbb_1\dd\int_\rrd(D_\vp b^*_\vp(T^n(f_0)))\cdot\nabla(T^n(f_0))d\xi
	{+}\(\dd\frac12T^{n-1}(f_0)+\frac12f,T^n(f_0) \)_2 . \nonumber\end{eqnarray}
We set 
\begin{equation}\label{b2.27v}
	g_\vp(r):=\int^r_0 b^*_\vp(\tau)d\tau,\ r\in\rr.\end{equation}
By Hypothesis (iii) we have 
$0\le g_\vp(r)\le|b^*_\vp|_\9r,\ r\in\rr,$ and so the right hand side of \eqref{b2.27iv} is equal to

$$-\lbb_1({\rm div}\ D_\vp,g_\vp(T^n(f_0)))_2+\(\frac12\,(T^{n-1}(f_0)+f),T^n(f_0)\)_2,$$
 where we recall that $\divv\,D_\vp\in L^2$ by \eqref{e2.8a}.  
Thus, by \eqref{e2.15aa} we obtain
$$\barr{r}
\dd|T^{n}(f_0)|^2_2\le\lbb_1|b^*_\vp|_\9|({\rm div}\,D_\vp)^-|_2|T^n(f_0)|_2
+\frac12|T^n(f_0)|^2_2\vsp
 +\dd\frac14(|T^{n-1}(f_0)|^2_2
 +|f|^2_2),\earr$$
therefore, 
 $$\dd |T^{n}(f_0)|^2_2\le C_\vp+\frac23\,|T^{n-1}(f_0)|^2_2,$$
 where 
 $$  C_\vp:=\dd\frac{16}{3}\,
\lbb^2_1|b^*_\vp|^2_\9|({\rm div}\,D_\vp)^-|^2_2+\frac23\,|f|^2_2.$$
 Consequently, by iteration 
$$|T^{n}(f_0)|^2_2\le C_\vp 
\sum^{n}_{k=0} \(\frac23\)^k+\(\frac23\)^n|f_0|^2_2,\ n\in\nn.$$
Hence, we get  
$$|J^\vp_\lbb(f)|^2_2\le\liminf_{n\to\9}|T^n(f_0)|^2_2\le 3C_\vp<\9,$$so, \eqref{b2.27'''} holds for $\lbb=2\lbb_1,$ and finally we get  \eqref{b2.27'''} for all $\lbb>0.$

 Let $\lbb^\vp_0:=(|b|_\9(|({\rm div}\,D)^-|_\9+\vp|\eta'|_\9|D|_\9))^{-1}$   and $M^\vp_\lbb:=1+\frac\lbb{\lbb^\vp_0}$.    	 
Then, for $f\in L^1\cap L^\9$ and $y_\vp:=J^\vp_\lbb(f),$
$\lbb>0$, we~have\begin{equation}
	\label{e2.23a2}
	|y_\vp|_\9\le
	\(1+\frac\lbb{\lbb^\vp_0}\)
	|f|_\9,\ \ff\lbb\in(0,\lbb^\vp_0). 
\end{equation}
 To prove \eqref{e2.23a2}, we note that for $\lbb\le\lbb^\vp_0$, we get by \eqref{e2.7}  and \eqref{e2.8a} that
\begin{eqnarray}
&& \hspace*{-14mm}(y_\vp{-}M^\vp_\lbb|f|_\9) 
	{+}\lbb G_\vp(\beta_\vp(y_\vp){-}\beta_\vp(|f|_\9))
	{+}\lbb\,\divv(D_\vp(b^*_\vp(y_\vp){-}b^*_\vp(|f|_\9)))\label{2.36'}\\[1mm]
&& \hspace*{8mm}
=f-M^\vp_\lbb|f|_\9-\lbb\psi(\vp)\beta_\vp(|f|_\9)-|f|_\9\lbb b_\vp(|f|_\9)\divv\,D_\vp\nonumber\\[1mm]
&& \hspace*{8mm}
\le f-|f|_\9(M^\vp_\lbb-\lbb|b|_\9(\divv\,D_\vp)^-)\nonumber\\[1mm]
&& \hspace*{8mm}
\le f-|f|_\9\(M^\vp_\lbb-\dd\frac\lbb{\lbb^\vp_0}\)
=f-|f|_\9\le0\nonumber\end{eqnarray}
because $G_\vp 1=\Psi(\vp)$, since   $\calf(1)=(2\pi)^{\frac d2}\,\delta_0$. 
 Taking scalar product in $L^2$ with $\calx_\delta((\beta_\vp(y_\vp)-\beta_\vp(|f|_\9)^+))$ in \eqref{2.36'}, letting $\delta\to0$ and using \eqref{e2.15a}  and a similar argument as that leading to \eqref{e2.22prim},  we get  
$$y_\vp\le \(1+\dd\frac\lbb{\lbb^\vp_0}\)
|f|_\9,\mbox{ a.e. in }\rrd,$$and, similarly, for $-y_\vp$ which yields \eqref{e2.23a2} for $\lbb\in(0,\lbb^\vp_0)$.  So, altogether we have 
\begin{equation}
	\label{e2.23aa}
	|J^\vp_\lbb(f)|_1+|J^\vp_\lbb(f)|_\9\le |f_1|+2|f|_\9,\ \ff\vp>0, \ \lbb\in(0,\lbb_0). 
\end{equation}

Now, fix $\lbb\in(0,\lbb^1_0)$ and $f\in L^1\cap L^\9$. As above, for $\vp\in(0,1]$ set
$y_\vp:=J^\vp_\lbb(f).$

Then, since $\beta_\vp(y_\vp)\in H^1$ and $G_\vp\beta_\vp(y_\vp)\in L^2$ by \eqref{b2.27'}, 
\begin{equation}\label{2.38}
	\barr{l}
		(y_\vp,\beta_\vp(y_\vp))_2+
	\lbb(G_\vp(\beta_\vp(y_\vp)),\beta_\vp(y_\vp))_2\vsp\qquad 
	=\lbb\dd\int_\rrd(D_\vp b^*_\vp(y_\vp))\cdot\nabla\beta_\vp(y_\vp)dr+(f,\beta_\vp(y_\vp))_2.\earr 
\end{equation}
Setting
\begin{equation}\label{2.39}
	\wt\Psi_\vp(r):=\int^r_0 b^*_\vp(\tau)\beta'_\vp(\tau)d\tau,\ r\in\rr,
\end{equation}by Hypotheses (iii), (iv)  we have
$$0\le\wt\Psi_\vp(r)\le\frac12\,|b|_\9(|\beta'|_\9+1)r^2,\ \ff r\in\rr,$$and hence, since $y_\vp\in H^1$, the right hand side of \eqref{2.38} is equal to
$$-\lbb\int_\rrd{\rm div}\ D_\vp\wt\Psi_\vp(y_\vp)dx+(f,\beta_\vp(y_\vp))_2,$$which, because $(y_\vp,\beta_\vp(y_\vp))_2\ge0$   and $H^1\subset H^\Psi\subset H^s$ (see \eqref{2.0}), by \eqref{e2.8a} and Hypothesis (iv) implies that 
$$\lbb|\beta_\vp(y_\vp)|^2_{H^{\Psi,\vp}}
\le\frac12\,\lbb|b|_\9(|\beta'|_\9+1)
(|({\rm div}\ D)^-|_\9+|\eta'|_\9|D|_\9)|y_\vp|^2_2
+\frac12\,|\beta_\vp(y_\vp)|^2_2+\frac12\,|f|^2_2.$$
Since $|\beta_\vp(r)|\le({\rm Lip}(\beta)+1)|r|,\ r\in\rr,$ by \eqref{e2.23aa} we obtain
\begin{equation}\label{2.40}
\sup_{\vp\in(0,1]}|\beta_\vp(y_\vp)|^2_{H^{\Psi,\vp}}\le C\ \mbox{ for some $C\in(0,\9)$.}\end{equation}
Obviously, we have for all $u\in H^\Psi$ $(\subset\dot H^\Psi)$ and all $\vp\in(0,1]$
\begin{equation}\label{2.41}
	|\Psi(-\Delta)^{\frac12}u|^2_2\le|\Psi(\vp I-\Delta)^{\frac12}u|^2_2\le|\Psi(-\Delta)
	^{\frac12}u|^2_2+\Psi(\vp)|u|^2_2,
\end{equation}where we used the sub linearity of $\Psi$ in the second step.

Hence, we conclude from  \eqref{e2.23aa}  and \eqref{2.40}  that  (along a subsequence) as $\vp\to0$,
\begin{eqnarray*}
\beta_\vp(y_\vp)&\to& z\mbox{ weakly in }H^\Psi,  \mbox{ hence strongly in $L^2_{\rm loc}(\rrd)$ (by \eqref{2.0})},\\
G_\vp(\beta_\vp(y_\vp))&\to&\Psi(-\Delta)z\mbox{ in }S',\\
y_\vp&\to&y\mbox{ weakly in $L^2$ and weakly$^*$ in $L^\9$},\end{eqnarray*}
(where the second statement follows, because 
$G_\vp(\vf)\to\Psi(-\Delta)\vf$ in $L^2$ for all $\vf\in S.$ Hence (selecting another subsequence, if necessary),
$\beta(y_\vp)\to z,$ a.e.  
Since $\beta\1$ is continuous, it follows that 
$y_\vp\to\beta\1(z)=y,\mbox{ a.e.},$  and, therefore, $z=\beta(y)\in H^\Psi$.   Furthermore, we have
$$b^*_\vp(y_\vp)\to b^*(y)\mbox{ weakly in }L^2.$$Recalling that $y_\vp$ solves \eqref{e2.7}, we can let $\vp\to0$ in \eqref{e2.7} to find that 
\begin{equation}\label{2.27'''}
	y+\lbb\Psi(-\Delta)\beta(y)+\lbb\,{\rm div}(Db^*(y))=f\mbox{ in }S'.\end{equation}Since $\beta\in {\rm Lip}(\rr)$, the operator $(A_0,D(A_0))$ defined in \eqref{e1.5}  is closed in $L^1$. If $y$ is as in \eqref{2.27'''}, we define
$$J_\lbb(f):=y\in D(A_0),\ \lbb\in(0,\lbb_0).$$  
Moreover,  by  Lemma A in the Appendix (by choosing $\lbb_0$ a little smaller, if necessary),      \eqref{2.27'''}  has a unique solution $y\in D(A_0)\cap L^1\cap L^\9$ and so, 
\begin{equation}\label{e2.40prim}
	\lim_{\vp\to0}\,y_\vp=y\mbox{ weakly in $L^2$ and weakly$^*$ in $L^\9$}\end{equation}
	(not only along a subsequence, which might depend on $f$ and $\lbb$) and, for every sequence $\vp_n\to0$,
	\begin{equation}\label{e2.40secund}
	\lim_{k\to\9}\,y_{\vp_{n_k}}=y,\mbox{ a.e.  for some subsequence }\{\vp_{n_k}\}.
	\end{equation}  
Then \eqref{e2.6a} holds and it follows by \eqref{e2.19} (which, as mentioned earlier, holds, in fact, for all $\lbb>0$), \eqref{e2.40secund} and Fatou's lemma that for $f_1,f_2\in L^1\cap L^\9$
\begin{equation}\label{2.27iv}
	|J_\lbb(f_1)-J_\lbb(f_2)|_1\le|f_1-f_2|_1.\end{equation}
Hence $J_\lbb$ extends continuously to all of $L^1$, still satisfying \eqref{2.27iv} for all $f_1,f_2\in L^1$. Then it follows by the closedness of $(A_0,D(A_0))$ on $L^1$ that $J_\lbb(f)\in D(A_0)$ and that it solves \eqref{2.27'''} for all $f\in L^1$.  So, properties \eqref{e2.1}, \eqref{e2.2} and the last statement in Lemma \ref{l2.1} are proven. 

Clearly, \eqref{e2.5} and \eqref{e2.5a} follow from \eqref{e2.13} and \eqref{e2.23a2}, respectively. 

Hence, Lemma \ref{l2.1} is proved except for \eqref{e2.3} and \eqref{e2.4}.  So, let $f\in L^1\cap L^\9;\ \lbb_1,\lbb_2\in(0,\lbb_0)$. Then 
$$(I+\lbb_1A_0)J_{\lbb_2}(f)=\frac{\lbb_1}{\lbb_2}\,f+\(1-\frac{\lbb_1}{\lbb_2}\)J_{\lbb_2}f.$$
 By the final statement of Lemma \ref{l2.1}, $J_{\lbb_1}\(\frac{\lbb_1}{\lbb_2}\,f+\(1-\frac{\lbb_1}{\lbb_2}\)J_{\lbb_2}\,f\)$ is the only solution to this equation in $D(A_0)\cap L^1\cap L^2$. But also $J_{\lbb_2}\,f$ is in the latter space, so \eqref{e2.3} follows. 
 
Now, let us prove \eqref{e2.4}. We may assume that $f\in L^1\cap L^\9$ and set $y:=J_\lbb(f)$. Let $\psi_n\in C^\9_0(\rrd),\ \psi_n\uparrow1$, as $n\to\9$, $\lim\limits_{n\to\9}\nabla \psi_n=0$ on $\rrd$, with $\sup\limits_n|\nabla\psi_n|_\9<\9$. Define
$$\vf_n:=(I+\Psi(-\Delta))\1\psi_n=g^{\Psi}_1*\psi_n,\ n\in\nn,$$ where $g^\Psi_\vp$ is as in the 
Appendix. Then, we have by (A.2) 
\begin{eqnarray}
	\label{b2.36}
	&\vf_n\uparrow1,\ \nabla\vf_n\to0\mbox{ on }\rrd\ \mbox{ as }n\to\9,\\
	&\dd\sup\limits_n(|\vf_n|_\9+|\nabla\vf_n|_\9)<\9,\ \
	\vf_n\in L^1\cap H^{2s},\ n\in\nn.\nonumber
\end{eqnarray}
Furthermore,
$$\Psi(-\Delta)\vf_n=\psi_n-(I+\Psi(-\Delta))\1\psi_n\in L^1\cap L^\9,$$are bounded in $L^\9$ and, as $n\to\9$,
$\Psi(-\Delta)\vf_n\to0\ dx-\mbox{a.e.}$ Hence, 
\begin{equation}\label{b2.37}
	\lim_{n\to\9}\int_\rrd\Psi(-\Delta)\vf_n\,\beta(y)dx=0.\end{equation}
Consequently, since $\beta(y)\in L^1$, $y\in D(A_0)$ with $A_0y\in L^1$, we have
$$\barr{l}
 \dd\int_\rrd A_0y\,dx
=\lim_{n\to\9}\int_\rrd\vf_n A_0y\,dx 
=-\int_\rrd\beta(y)dx\vsp
\qquad+\dd\lim_{n\to\9}
\,{}_{S}\!\<\vf_n,(I+\Psi(-\Delta))\beta(y)+{\rm div}(Db^*(y))\>_{S'}\\
\qquad=-\dd\int_\rrd\!\beta(y)dx{+}\lim_{n\to\9}
\int_\rrd(I+\Psi(-\Delta))\vf_n\,\beta(y)dx\vsp
\qquad  
=+\dd\lim_{n\to\9}\int_\rrd\!\!\nabla\vf_n\cdot Db^*(y)dx,
\end{array}$$which by \eqref{b2.36} and \eqref{b2.37} is equal to zero. Hence, integrating the equation
$y+\lbb A_0y=f$ over $\rrd$, \eqref{e2.4} follows.
\end{proof}

Now, define the operator $A$ by
\begin{eqnarray}\label{2.35}
&	D(A) :=J_\lbb(L^1)\ (\subset D(A_0)),\  \lbb\in(0,\lbb_0);\
	A(y) :=A_0(y),\ y\in D(A). \end{eqnarray}
Then, by \eqref{e2.3},   $J_\lbb(L^1)$ is independent of $\lbb\in(0,\lbb_0)$ and,  clearly,  
 $J_\lbb=(I+\lbb A)\1,\ \lbb\in(0,\lbb_0).$ 
Therefore, we have
\begin{lemma}\label{l2.2} Under Hypotheses  {\rm(i)--(iv)}, the operator $A$ is $m$-accretive in $L^1$ and $(I+\lbb A)\1=J_\lbb$, $\lbb\in(0,\lbb_0)$. Moreover, if $\beta\in C^\9(\rr)$, then $\ov{D(A)}=L^1$.\end{lemma}
\n(Here, $\ov{D(A)}$ is the closure of $D(A)$ in $L^1$.) 

 \begin{proof}  It only remains to prove the last statement. So, let   $\beta\in C^\9(\rr)$.Then, by assumption (ii),  
$$A_0(\vf)=\Psi(-\Delta)\beta(\vf)+{\rm div}(Db(\vf)\vf)\in L^1,\ \ff\vf\in C^\9_0(\rrd),$$and so $\ov{D(A)}=L^1$, as claimed.\end{proof}
 
Then, by  the Crandall \& Liggett theorem (see, e.g., \cite{1}, p.~131), the Cauchy problem \eqref{e1.4} has, for each $u_0\in\ov{D(A)}$, a unique mild solution $u=u(t,u_0)\in C([0,\9);L^1)$ and $S(t)u_0=u(t,u_0)$ is a $C_0$-semigroup of contractions on $L^1$. 
Moreover, by \eqref{e2.5a} and the exponential formula
$$S(t)u_0=\lim_{n\to\9} \(J_{\frac tn}\)^{n}u_0,\ \ff\,t\ge0,$$
it follows that 
$$|S(t)u_0|_\9\le
e^{|(\divv\,D)^-|_\9|b|_\9t} |u_0|_\9,\ \ff t\ge0.$$ Hence, if $u_0\in L^1\cap L^\9,$ then 
$S(t)u_0\in L^\9((0,T)\times\rrd)$, $T>0$. Furthermore, by \eqref{e2.5} if $u_0\ge0$, then $S(t)u_0\ge0$. 

Let us now show   that $u=S(t)u_0$ is a Schwartz distributional solution, that is, \eqref{e1.10} holds. By \eqref{e16}, we have
\begin{eqnarray*}
  \int^\9_0\!\!\!dt\Bigg(\int_\rrd\vf(t,x)\frac1h(u_h(t,x)-u_h(t-h,x))dx  + \!\!\int_\rrd
(\vf(t,x)\Psi(-\Delta)\beta(u_h(t,x))\\
 -\nabla_x\vf(t,x)\cdot D(x)b^*(u_h(t,x)))dx\Bigg)=0,\  \ff\vf\in C^\9_0([0,\9)\times\rrd).\end{eqnarray*}
This yields
\begin{eqnarray*}
&& \int^\9_0dt
\Bigg(\int_\rrd u_h(t,x)\frac1h(\vf(t+h,x)-\vf(t,x))dx\\
&&\quad+\dd\int_\rrd(\beta(u_h(t,x))\Psi(-\Delta)\vf(t,x)
-\nabla_x\vf(t,x)\cdot D(x)b^*(u_h(t,x)))dx\Bigg)\\
&&\quad+\dd\frac1h\int^h_0dt\int_\rrd u_0(x)\vf(t,x)dx=0,\ \ff\vf\in C^\9_0([0,\9)\times\rrd).
\end{eqnarray*}Taking into account that, by \eqref{e1.6} and assumptions (i)--(iii), $\beta(u_h)\to\beta(u)$, $b^*(u_h)\to b^*(u)$ in $C([0,T];L^1)$ as $h\to0$ for each $t>0$, we get  that \eqref{e1.10} holds. 

We have, therefore, the following existence result for equation~\eqref{e1.1}.
\begin{theorem}\label{t2.3} Assume $s\in\(\frac12,1\)$ and that Hypotheses {\rm(i)--(v)} hold. Then, there is a $C_0$-semigroup of contractions $S(t):L^1\to L^1$, $t\ge0$, such that for each \mbox{$u_0\in \ov{D(A)}$,} which is $L^1$ if $\beta\in C^\9(\rr)$, $u(t,u_0)=S(t)u_0$ is a mild  solution to \eqref{e1.1}.  Moreover, if $u_0\ge0$, a.e. on $\rrd$, 	
	\begin{eqnarray}
		& 	u(t,u_0)\ge0,\mbox{\ \ a.e. on }\rrd,\ \ff\,t\ge0,\label{e2.28}\\ 
		&\dd\int_\rrd u(t,u_0)(x)dx=\int_\rrd u_0(x)dx,\ \ \ff\,t\ge0.
		\label{e2.29}\end{eqnarray} 	
	Moreover, $u$ is a distributional solution to \eqref{e1.1} on $[0,\9)\times\rrd.$ Finally, if $u_0\in L^1\cap L^\9$,  then all above assertions remain true, if we drop the assumption $\beta\in{\rm Lip}(\rr)$ from Hypothesis {\rm(i)}, and additionally we have that \mbox{$u\in L^\9((0,T)\times\rrd)$,} $T>0$.\end{theorem}

\section{The uniqueness of distributional solutions}\label{s3}
\setcounter{equation}{0}

In general, the mild solution $u$ given by Theorem \ref{t2.3} is not unique because the operator $A$ constructed in Lemma \ref{l2.2}  depends on the approximating operators $A_\vp y\equiv\Psi(\vp I-\Delta)\beta_\vp(y)+{\rm div}(D_\vp b_\vp(y)y)$ and so $u=S(t)u_0$ may be viewed as a {\it viscosity-mild} solution to \eqref{e1.1}. However, as seen here, this mild solution is even unique in the much larger class of Schwartz distributional solutions under the following hypotheses on $\beta,b$ and $D$: 
\begin{itemize}
\item[\rm(j)] $\beta\in C^1(\rr),\ \beta'(r)>0,\ \ff\,r\in\rr,\ \beta(0)=0.$\vspace*{-2mm}
\item[\rm(jj)] $D\in L^\9(\rrd;\rrd).$\vspace*{-2mm}
\item[\rm(jjj)] $b\in C^1(\rr).$
\end{itemize}

\begin{theorem}\label{t3.1} Assume that $0<T<\9$, $d\ge2,$ and that  Hypotheses {\rm(j)--(jjj)} and {\rm(v)} hold.  Let $y_1,y_2\in L^\9((0,T){\times}\rrd)$ be two distributional solutions to \eqref{e1.1} on $(0,T)\times\rr^d$ $($in the sense of \eqref{e1.10}$)$   such that  $y_1{-}y_2\in L^1((0,T)\times\rrd)$ $\cap L^\9(0,T;L^2)$ and 
	\begin{equation}
		\label{e3.1}\lim_{t\to0}\ {\rm ess} \sup_{\hspace*{-4mm}s\in(0,t)}|(y_1(s)-y_2(s),\vf)_2|=0,\ \ff\vf\in C^\9_0(\rrd).
	\end{equation}	Then $y_1\equiv y_2$. If $D\equiv0$, then Hypothesis {\rm(j)} can be relaxed to
	\begin{itemize}
		\item[\rm(j)$'$] $\beta\in C^1(\rr),\ \beta'(r)\ge0,\ \ff\,r\in\rr,\ \beta(0)=0.$
\end{itemize}\end{theorem} 

\begin{proof} (The proof is similar to that of  Theorem 3.1 in \cite{a9}, but it has to be adapted substantially.) Replacing, if necessary, the functions $\beta$ and $b$ by
$$\beta_N(r)=\left\{\barr{ll}
\beta(r)&\mbox{ if }|r|\le N,\vsp
\beta'(N)(r-N)+\beta(N)&\mbox{ if }r>N,\vsp
\beta'(-N)(r+N)+\beta(-N)&\mbox{ if }r<-N,\earr\right.$$ 
$$ b_N(r)=\left\{\barr{ll}
b(r)&\mbox{ if }|r|\le N,\vsp
b'(N)(r-N)+b(N)&\mbox{ if }r>N,\vsp
b'(-N)(r+N)+b(-N)&\mbox{ if }r<-N,\earr\right.$$where $N\ge\max\{|y_1|_\9,|y_2|_\9\}$, by (j) we may assume that 

\begin{equation}\label{e3.2}
	\beta',b'\in C_b(\rr),\ \beta'\ge\alpha_2\in(0,\9)\end{equation}and, therefore,   we have
\begin{eqnarray}
	\alpha_1|\beta(r)-\beta(\bar r)|&\ge& 
	|b^*(r)-b^*(\bar r)|,\ \ \ff\,r,\bar r\in\rr,\label{e3.3}\\
	(\beta(r)-\beta(\bar r))(r-\bar r)&\ge&\alpha_3|\beta(r)-\beta(\bar r)|^2,\   \ff\,r,\bar r\in\rr,\label{e3.3a}\end{eqnarray}
where $b^*(r)=b(r)r$, $\alpha_1\ge0,$ and $\alpha_3:=|\beta'|^{-1}_\9$. We set
\begin{eqnarray}
	&	\Phi_\vp(y)=(\vp I+\Psi(-\Delta))\1y,\ \ff\,y\in L^2,\label{e3.4}\\
& z=y_1-y_2,\ w=\beta(y_1)-\beta(y_2),\ b^*(y_i)\equiv b(y_i)y_i,\ i=1,2.\nonumber\end{eqnarray}
We have that $\Phi_\vp:L^p\to L^p$, $\ff p\in[1,\9]$ and
\begin{equation}\label{e3.5}
	\vp|\Phi_\vp(y)|_p\le |y|_p,\ \ \ff y\in L^p,\ \vp>0.
\end{equation}
 Furthermore,
$$z_t+\Psi(-\Delta) w+{\rm div}\,  D(b^*(y_1)-b^*(y_2))=0\mbox{ in }\cald'((0,T)\times\rrd).$$
We set
\begin{equation}\label{e3.6a}
	z_\vp=z*\theta_\vp,\ w_\vp=w*\theta_\vp,\ \zeta_\vp=(D(b^*(y_1)-b^*(y_2)))*\theta_\vp,\end{equation}where $\theta\in C^\9_0(\rr^d),$ $\theta_\vp(x)\equiv\vp^{-d}\theta\(\frac x\vp\)$ is a standard mollifier. 
	
	We note that $z_\vp,w_\vp,\zeta_\vp,\Psi(-\Delta)w_\vp,{\rm div}\,\zeta_\vp\in L^2(0,T;L^2)$ and  we have
\begin{equation}\label{e3.6}
	(z_\vp)_t+\Psi(-\Delta) w_\vp+{\rm div}\,\zeta_\vp\ \mbox{ in }\cald'(0,T;L^2) .\end{equation}
This yields  $\Phi_\vp(z_\vp),\Phi_\vp(w_\vp),{\rm div}\,\Phi_\vp(\zeta_\vp)\in L^2(0,T;L^2)$ and 
\begin{equation}\label{e3.7}
	\barr{ll}
	(\Phi_\vp(z_\vp))_t=-\Psi(-\Delta)\Phi_\vp(w_\vp)-{\rm div}\Phi_\vp(\zeta_\vp)\mbox{ in }\cald'(0,T;L^2).\earr\end{equation}
By \eqref{e3.6}, \eqref{e3.7} it follows that $(z_\vp)_t=\frac d{dt}\,z_\vp$,   $(\Phi_\vp(z_\vp))_t=\frac d{dt}\,\Phi_\vp(z_\vp)\in L^2(0,T;L^2)$. 
This implies that $z_\vp, \Phi_\vp(z_\vp)\in H^1(0,T;L^2)$ and both $[0,T]\ni t\mapsto z_\vp(t)\in L^2$ and $[0,T]\ni t\mapsto\Phi_\vp(z_\vp(t))\in L^2$ are absolutely continuous. Moreover, it follows by \eqref{e3.5} and \eqref{e3.7}  that
\begin{equation}\label{e3.8}
	\Phi_\vp(z_\vp),\Phi_\vp(w_\vp)\in L^2(0,T; L^2).\end{equation}
We set $h_\vp(t)=(\Phi_\vp(z_\vp(t)),z_\vp(t))_2$ and get 
\begin{equation}
\hspace*{-4mm}\barr{ll}
	h'_\vp(t)\!\!\! &=2(z_\vp(t),(\Phi_\vp(z_\vp(t)))_t)_2\\
	&=2(\vp\Phi_\vp(w_\vp(t)){-}w_\vp(t){-}{\rm div}\Phi_\vp(\zeta_\vp(t)),z_\vp(t))_2\label{e3.9}\\
&=2\vp(\Phi_\vp(z_\vp(t)),w_\vp(t))_2{+}2(\nabla\Phi_\vp(z_\vp(t)),\zeta_\vp(t))_2 
 -2(z_\vp(t),w_\vp(t))_2,\\
 &\hfill\mbox{ a.e. }t\in(0,T).\earr
\end{equation}
By \eqref{e3.7}--\eqref{e3.9} it follows that \mbox{$t\to h_\vp(t)$} has an absolutely continuous $dt$-version on $[0,T]$ which we shall consider from now on. Since, by \eqref{e3.3a}, we have 
\begin{equation}\label{e3.10}
(z_\vp(t),w_\vp(t))_2\ge\alpha_3|w(t)|^2_2+\gamma_\vp(t),\end{equation}
where
\begin{equation}\label{e3.11}
\gamma_\vp(t):=(z_\vp(t),w_\vp(t))_2-(z(t),w(t))_2,\end{equation} 
we get, by \eqref{e3.3} and \eqref{e3.8}, \newpage
\begin{eqnarray}\label{e3.12}
		0&\le& h_\vp(t)
	\le h_\vp(0+)+2\vp\!\dd\int^t_0 
	(\Phi_\vp(z_\vp(s)),w_\vp(s))_2ds\\[-2mm]
&	-&2\alpha_3\!\dd\int^t_0
	|w(s)|^2_2ds 
 +2\alpha_1|D|_\9\dd\int^t_0\!|\nabla\Phi_\vp(z_\vp(s))|_2|w(s)|_2ds\nonumber\\[-2mm]
	&+&2
	\int^t_0|\gamma_\vp(s)|ds,\, \ff t\in[0,T].\nonumber\end{eqnarray} Moreover, since $z\in L^\9((0,T)\times\rrd)$,   by \eqref{e3.5} we have
\begin{equation}
	\label{e3.13}
	\vp|\Phi_\vp(z_\vp(t))|_\9\le|z_\vp(t)|_\9\le|z(t)|_\9,\mbox{\ \ a.e. }t\in(0,T).\end{equation} 
As $t\to\Phi_\vp(z_\vp(t))$ has an $L^2$ conti\-nuous version on $[0,T]$, there exists $f\in L^2$ such that 
 $\lim\limits_{t\to0}\Phi_\vp(z_\vp(t))=f\mbox{\ \ in }L^2.$ 
Furthermore, for every $\vf\in C^\9_0(\rrd)$, $s\in(0,T),$
$$0\le h_\vp(s)\le|\Phi_\vp(z_\vp(s))-f|_2|z_\vp(s)|_2+|f-\vf|_2|z_\vp(s)|_2+|(\vf*\theta_\vp,z(s))_2|,$$ and so, by \eqref{e3.1},
\begin{eqnarray*}
0\le h_\vp(0+) 
&=&\dd\lim_{t\downarrow0}h_\vp(t)
=\lim_{t\to0}\
{\rm ess}\sup_{\hspace*{-4mm}s\in(0,t)}h_\vp(s)\\[-2mm]
&\le&\left(\dd\lim_{t\to0}|\Phi_\vp(z_\vp(t))-f|_2+|f-\vf|_2\right)|z_\vp|_{L^\9(0,T;L^2)}\\
&&+\lim_{t\to0}\ {\rm ess}\sup_{\hspace*{-4mm}s\in(0,t)}|(\vf*\theta_\vp,z(s))_2|=|f-\vf|_2|z_\vp|_{L^\9(0,T;L^2)}.\end{eqnarray*}
Since $C^\9_0(\rrd)$ is dense in $L^2(\rrd)$, we find that
\begin{equation}\label{c3.16'}
	h_\vp(0+)=0.\end{equation}
On the other hand, taking into account that, for a.e. $t\in(0,T)$,
\begin{equation}\label{e3.15}
	\vp\Phi_\vp(z_\vp(t))+\Psi(-\Delta)\Phi_\vp(z_\vp(t))=z_\vp(t),\end{equation}we get that, for a.e. $t\in(0,T)$,
\begin{eqnarray}
	&\vp|\Phi_\vp(z_\vp(t))|^2_2+|(\Psi(-\Delta))^{\frac 12}\Phi_\vp(z_\vp(t))|^2_2=(z_\vp(t),\Phi_\vp(z_\vp(t)))_2=h_\vp(t),\quad\label{e3.16}\\
	&	\vp|(\Phi_\vp(z_\vp(t),w_\vp(t)))_2|\le\vp|\Phi_\vp(z_\vp(t))|_\9|w_\vp(t)|_1
	\le|z(t)|_\9|w(t)|_1.\quad\label{e3.17a}\end{eqnarray}
By \eqref{e3.15}, we have
\begin{equation}\label{e3.17aa}\calf(\Phi_\vp(z_\vp(t)))=(\vp+\Psi(|\xi|^2))^{-1}\calf(z_\vp(t)).\end{equation}
Therefore, by Parseval's formula,
\begin{eqnarray*}
 |\nabla\Phi_\vp(z_\vp(t))|^2_2&=& \int_\rrd
\frac{|\calf(z_\vp(t))(\xi)|^2|\xi|^2}
{(\vp+\Psi(|\xi|^2))^2}\ d\xi,\ \ff t\in(0,T),\\
  h_\vp(t)&=& \int_\rrd\ \frac{|\calf(z_\vp(t))(\xi)|^2}{\vp+\Psi(|\xi|^2)}\ d\xi,\ \ff t\in(0,T).\end{eqnarray*}
Then, by \eqref{e1.4aa} this yields for some $C\in(0,\9)$ independent of $\vp$
\begin{eqnarray}
 	|\nabla\Phi_\vp(z_\vp(t))|^2_2 
	&\le&	CR^{2(1-s)} 
	\int_{[|\xi|\le R]}\frac{|\calf(z_\vp(t))(\xi)|^2}{\vp+\Psi(|\xi|^2)}\ d\xi\label{e3.15z}\\ 
&& +   C \int_{[|\xi|\ge R]}|\calf(z_\vp(t))(\xi)|^2|\xi|^{2(1-2s)}d\xi \nonumber\\
	&\le& CR^{2(1-s)}h_\vp(t)+CR^{2(1-2s)}|z_\vp(t)|^2_2, \nonumber\end{eqnarray}
$\ff t\in(0,T),\ R>0,$ because $2s\ge1.$ 

We shall now prove that
\begin{equation}\label{e3.15a}
	\lim_{\vp\to0}\vp(\Phi_\vp(z_\vp(t)),w_\vp(t))_2=0,\mbox{ a.e. }t\in(0,T).\end{equation}
By \eqref{e3.17aa} and \eqref{e1.4aa} it follows  for some $C\in(0,\9)$ independent of $\vp$
\begin{eqnarray*}
 	|(\Phi_\vp(z_\vp(t)),w_\vp(t))_2| &=&|(\calf(\Phi_\vp(z_\vp(t)), \ov\calf(w_\vp(t)))_2|\\[1mm] 
&\le & C\dd\int_\rrd\frac{|\calf(z_\vp(t))(\xi)|\,
		|\calf(w_\vp(t))(\xi)|}
	{\vp+|\xi|^{2s}}\, d\xi\\ 
	& \le & C\(\dd\int_\rrd\left|\frac{\calf(z_\vp(t))}{\vp+|\xi|^{2s}}\right|^2d\xi\)^{\frac12}|w_\vp(t)|_2,\end{eqnarray*}
and since
$$\frac{\calf(z_\vp(t))}{\vp+|\xi|^{2s}}=\calf((\vp I+(-\Delta)^s)\1z_\vp(t)),\ t\in(0,T),$$
this yields 
\begin{equation}\label{3.20}
	|(\Phi_\vp(z_\vp(t)),w_\vp(t))_2|\le C|(\vp I+(-\Delta)^s)\1z_\vp(t)|_2|w(t)|_2. \end{equation}On the other hand, for each $f\in L^2(\rrd)$, by \cite[Appendix]{a9} we have
\begin{equation}\label{3.21}
	(\vp I+(-\Delta)^s)\1f(x)=\int_\rrd g^s_\vp(x-\xi)f(\xi)d\xi,\end{equation}
	where
$$g^s_\vp(x)=\int^\9_0e^{-\vp\tau}d\tau\int^\9_0\frac{e^{-\frac{|x|^2}{4r}}}{(4\pi r)^{\frac d2}}\ \eta^s_\tau(dr),$$
and $(\eta^s_\tau)_{\tau\ge0}$ is the one-sided stable semigroup of order $s\in\(\frac12,1\).$ 

By (A.4), (A.7) and (A.10) in \cite{a9}, we have
\begin{eqnarray}
	&\vp\dd\int_\rrd g^s_\vp(x)dx=1,\label{e3.22a}\\ 
	&g^s_\vp\in L^\9(\rrd\setminus B_R(0)),\ \ff R>0,\ \
	g^s_\vp(x)=\vp^{\frac{d-2s}{2s}}\ g^s_1(\vp^{\frac1{2s}}\,x),\label{e3.22aa}  \end{eqnarray}
where $B_R(0)$ is the ball of radius $R$ around the origin in $\rrd$.

Then, by \eqref{3.21}--\eqref{e3.22aa}, we have via the Young inequality
\begin{eqnarray}
	&&\hspace*{-5mm}\vp|(\vp I+(-\Delta)^s)\1z_\vp(t)|_2=\vp|g^s_\vp*z_\vp(t)|_2\label{3.24}\le\vp|g^s_\vp*z_\vp(t)|^{\frac12}_\9|g^s_\vp*z_\vp(t)|^{\frac12}_1 \qquad\\
	&&\qquad\le\vp^{\frac12}|z_\vp(t)|^{\frac12}_1
	\dd\sup_{x\in\rrd}
	\(\int_\rrd g^s_\vp(x-\xi)
	|z_\vp(t)(\xi)|d\xi\)^{\frac12}\nonumber\\
	&&\qquad\le\vp^{\frac d{4s}}\dd\sup_{x\in\rrd}\(
	\dd\int_\rrd g^s_1
	\(\vp^{\frac1{2s}}(x-\xi)\)|z_\vp(t)(\xi)|d\xi\)^{\frac12}|z(t)|^{\frac12}_1\le C_\delta\,\vp^{\frac d{4s}}|z(t)|_1 \nonumber\\
	&&\qquad+ \vp^{\frac d{4s}}|z(t)|^{\frac12}_\9
	\(\dd\int_{[\vp^{\frac1{2s}}|x-\xi|\le\delta]} g^s_1\(\vp^{\frac1{2s}}(x-\xi)\)d\xi\)^{\frac12}
	|z(t)|^{\frac12}_1,\nonumber\end{eqnarray}
where $C_\delta=\sup\{g^s_1(\xi);\,|\xi|> \delta\}$. Now, letting first $\vp\to0$ and then $\delta\to0$, \eqref{e3.15a} follows by \eqref{3.20} and \eqref{e3.22a}. 
Then, by \eqref{e3.17a} it follows that
\begin{equation}\label{e3.21}
	\lim_{\vp\to0}\vp\int^t_0(\Phi_\vp(z_\vp(s)),w_\vp(s))_2ds=0,\ t\in[0,T].\end{equation}
Next, by \eqref{e3.12}, \eqref{c3.16'} and \eqref{e3.15z}, we have  
\begin{eqnarray*}
0\le h_\vp(t)&\le& 2\vp\dd\int^t_0|(\Phi_\vp(z_\vp(r)),w_\vp(r))_2|dr-2\alpha_3\dd\int^t_0|w(r)|^2_2dr\\[-2mm]
&&+2\alpha_1|D|_\9\dd\int^t_0|\nabla\Phi_\vp(z_\vp(r))|_2|w(r)|_2dr 
+2\dd\int^t_0|\gamma_\vp(r)|dr\\[-2mm]
&\le&\eta_\vp(t)+2\alpha_1|D|_\9C^{\frac12}\dd\int^t_0
\(R^{1-s}h^{\frac12}_\vp(r)+R^{1-2s}|z_\vp(r)|_2\)
|w(r)|_2dr\\[-2mm]
&&-2\alpha_3\dd\int^t_0|w(r)|^2_2dr,\ \ff\, t\in[0,T],\ R>0,\end{eqnarray*}where
$$\eta_\vp(t):=2\vp\int^t_0|(\Phi_\vp(z_\vp(r)),w_\vp(r))_2|dr
+2\int^t_0|\gamma_\vp(r)|dr.$$
This yields
\begin{eqnarray}
&&	0 \le  h_\vp(t)
\le\eta_\vp(t){+}2\alpha_1|D|_\9C^{\frac12}
 \Big( R^{2(1-s)}\lbb \int^t_0  h_\vp(r)dr\label{e3.28a}\\ 
&&\qquad\qquad	+ \int^t_0 \Big(R^{1-2s}|z_\vp(r)|^2_2 
 +	   \Big(\frac1{4\lbb}+R^{1-2s}\Big)
	|w(r)|^2_2\Big)dr\Big)\nonumber\\ 	&&\qquad\qquad-2\alpha_3\dd\int^t_0|w(r)|^2_2dr,\ \ff\lbb>0,\ R>0. \nonumber\end{eqnarray}
Taking into account that, by \eqref{e3.2},  
\begin{equation}\label{e3.22'}
 |z_\vp(t)|_2\le 	|z(t)|_2\le \alpha^{-1}_2|w(t)|_2,\ \ \ff t\in(0,T),\end{equation}
we get for $\lbb$, $R>0$, large enough 
\begin{equation}\label{e3.22}
	\barr{r}
	\dd
	0\le h_\vp(t)\le\eta_\vp(t)+C\int^t_0 h_\vp(r)dr,\ \mbox{ for } t\in[0,T],\earr\end{equation}where $C>0$ is independent of $\vp$  and $\dd\lim_{\vp\to0}\eta_\vp(t)=0$ for all $t\in[0,T]$. 

By \eqref{e3.22}, it follows that
\begin{equation}\label{e3.23}
	0\le h_\vp(t)\le\eta_\vp(t)\exp(Ct),\ \ff\,t\in[0,T].\end{equation} 
This implies that $h_\vp(t)\to0$  as $\vp\to0$  for every $t\in[0,T]$, 
hence by \eqref{e3.16}   the left hand side of \eqref{e3.15} converges to zero in $S'$. Thus, $0=\lim\limits_{\vp\to0} z_\vp(t)=z(t)$ in $S'$ for a.e. $t\in(0,T)$, which implies $y_1\equiv y_2.$ If $D\equiv0$, we see by \eqref{e3.12} and  \eqref{c3.16'} that $0\le h_\vp(t)\le\eta_\vp(t)$, $\ff\,t\in(0,T)$, and so by \eqref{e3.21}  the conclusion follows without invoking that $\beta'>0$, which was only used to have 
\eqref{e3.22'}.
\end{proof} 


Similarly as Theorem \ref{t3.1}, one also obtains the li\-nearized uniqueness for equation \eqref{e1.10}.    

\begin{theorem}\label{t3.2} Under the assumptions of Theorem {\rm\ref{t3.1}}, let $T>0$, $u\in L^\9((0,T)\times\rrd)$  and let $y_1,y_2\in L^\9((0,T)\times\rrd)$ with $y_1{-}y_2\in L^1((0,T)\times\rrd)$ $\cap L^\9(0,T;L^2)$ be two distributional solutions to the equation	
	\begin{eqnarray}
	&&	\label{e3.30a} 
		y_t+\Psi(-\Delta)\(\dd\frac{\beta(u)}uy\)+{\rm div}(yDb(u))=0\mbox{ on }(0,T)\times\rrd,\\[-2mm]
	&&	y(0)=u_0,\nonumber
	\end{eqnarray}
	where $u_0$ is a measure of finite variation on $\rrd$ and $\frac{\beta(0)}0:=\beta'(0)$. If \eqref{e3.1} holds, then $y_1\equiv y_2$.\end{theorem}

\begin{proof} The proof is essentially the same as that of Theorem \ref{t3.1}   and, therefore, it will be sketched only. We note first that, by Hypotheses {\rm(j)--(jjj)},
$$\barr{c}
 \dd\frac{\beta(u)}u,\ b(u)\in L^\9((0,T)\times\rrd),\vsp
|Db(u)|_\9\le C_1\,\left|\dd\frac{\beta(u)}u\right|_\9,\  \left|\dd\frac{\beta(u)}u\right|\ge\alpha_2,\mbox{ a.e. in $(0,T)\times\rrd$}.\earr$$If $z=y_1-y_2,\ w=\frac{\beta(u)}u\,(y_1-y_2)$, we have therefore
\begin{eqnarray}
	wz&\ge&\alpha_2|w|^2, \mbox{ a.e. on $(0,T)\times\rrd,$}\label{e3.33}\\ 
	|Db(u)z|&\le& C_2|w|, \ \,\mbox{ a.e. on $(0,T)\times\rrd.$}\label{e3.34}
\end{eqnarray}
 We have 
$$z_t+\Psi(-\Delta)w+{\rm div}(Db(u)z)=0\mbox{ in }(0,T)\times\rrd,$$
and this yields (see \eqref{e3.6})
\begin{equation}\label{e3.38}
	(z_\vp)_t+\Psi(-\Delta)w_\vp+{\rm div}\,\zeta_\vp=0\mbox{ in }(0,T)\times\rrd,\end{equation} 
	\n where $z_\vp,w_\vp$ are as in \eqref{e3.6a} and $\zeta_\vp=(D(b(u))z)*\theta_\vp$.   If $\Phi_\vp$ is given by \eqref{e3.4}, by \eqref{e3.38} we get \eqref{e3.7} and, if $h_\vp(t)=(\Phi_\vp(z_\vp(t),z_\vp(t)))$, then   \eqref{e3.9} follows and so, by \eqref{e3.33} we get also in this case the estimates \eqref{e3.10} and \eqref{e3.12}.  From now on the proof is exactly the same as that of Theorem \ref{t3.1}. Namely, one gets that $h_\vp(0_+)=0$ and also that \eqref{e3.15z} and \eqref{e3.21} hold. Finally, one gets \eqref{e3.28a} and,  taking into account \eqref{e3.22'}--\eqref{e3.22}, one obtains that \eqref{e3.23} holds and so $z\equiv0$, as claimed.\end{proof}

\section{Applications to corresponding nonlinear\\ martingale problems}\label{s4}
\setcounter{equation}{0}

Here, we fix $T\in(0,\9)$ and need the following additional hypotheses on $\Psi$:
\begin{itemize}
	\item [\rm(vi)] $\dd\int^\9_1\log s\ \mu(ds)<\9.$
\end{itemize}
\n This is e.g. fulfilled if $\mu$ is as in \eqref{e1.7a}. 

\subsection{Existence}\label{s4.1}

Assume that Hypotheses (i)--(vi) hold and consider the nonlocal Kolmogorov operator corresponding to \eqref{e1.1}, i.e.,
\begin{equation}\label{e4.1}
	K_{u(t)} f(x)=\frac{\beta(u(t,x))}{u(t,x)}
	(\Psi(-\Delta)f)(x)+b(u(t,x))D(x)\cdot\nabla f(x),\  x\in\rrd,
\end{equation}
where $f\in C^2_c:=C^2_c(\rrd)$ and $u$ is the solution to \eqref{e1.1} from Theorem \ref{t2.3} with initial condition $u_0\in L^1\cap L^\9$.

By Theorem 13.6 in \cite{12''} it follows that, for all $f\in C^2_c$,
\begin{eqnarray}
	(\Psi(-\Delta)f)(x)\!\!\!
	&=&\int_{(0,\9)}\int_\rrd(f(x)-f(x+z))\frac1{\sqrt{4\pi t}}\,e^{-\frac{|z|^2}{4t}}\,dz\mu(dt)\nonumber\\
	&=&\int_\rrd(f(x)-f(x+z))\nu(|z|)dz,\label{e4.2}\\
\nu(r)&=&\int_{(0,\9)}\frac1{\sqrt{4\pi t}}\ e^{-\frac{r^2}{4t}}\,\mu(dt),\ r\in(0,\9).\nonumber\end{eqnarray}
Then, since  $\frac{\beta(u)}u\in L^\9((0,T)\times\rrd)$,  by (vi) it is easily seen that $K_{u(t)}$ is a Kolmogorov operator of the type considered in Section 1.2 in \cite{7aa}, which satisfies condition (1.18) in \cite{7aa}. Thus, by Theorem 1.5 ("{\it superposition principle}") in \cite{7aa}  and Remark 1.6 in \cite{7aa}, we get  
\begin{theorem}\label{t4.1}  Assume  Hypotheses {\rm(i)--(vi)}  hold and let $u(t,x),$ $t\in[0,T]$, $x\in\rrd$, be the solution from Theorem {\rm\ref{t2.3}}. Then, there exists a probability measure $\mathbb{P}$ on the Skorohod space $D([0,T];\rrd)$, which is a solution to the martingale problem corresponding to $(K_{u(t,\cdot)},C^2_c)$ in the sense of Definition {\rm1.3} in {\rm\cite{7aa}} with one dimensional time marginal densities given by $u(t,x),$ $t\in[0,T]$, $x\in\rrd$, i.e. for the canonical process $X_t$, $t\in[0,T]$, defined by $X_t(w)=w(t)$, $w\in D([0,T];\rrd)$, we have
	\begin{equation}\label{e4.3}
		(\mathbb{P}\circ X^{-1}_t)(dx)=u(t,x)dx,\ t\in[0,T].
	\end{equation}
\end{theorem}
\begin{definition}
	\label{d4.2} \rm The probability measure $\mathbb{P}$ in Theorem \ref{t4.1} is called a {\it solution to the nonlinear martingale problem corresponding to $(K_\Box,C^2_c)$, if} it solves the linear martingale problem corresponding to $(K_{\call_{X_t}},C^2_c)$  in the sense of  \cite[Definition 1.3]{7aa}, where $\call_{X_t}$ is its own one dimensional time marginal law density and $K_{\call_{X_t}}$ is defined as  in \eqref{e4.1} with $\call_{X_t}$ replacing $u(t,\cdot)$, $t\in[0,T]$. \end{definition}

\begin{remark}\label{r4.3}\rm We refer to the pioneering work \cite{a11'}, where such nonlinear martingale problems were studied for local Kolmogorov operators.
\end{remark}

\subsection{Uniqueness}\label{s4.2}

\begin{theorem}	\label{t4.3} Assume that Hypotheses {\rm(j)--(jjj), (v), (vi)} $($respectively, {\rm(j)$'$, (jj), (jjj), (v), (vi)} if $D\equiv0)$ hold. Let $\mathbb{P}$, $\wt{\mathbb{P}}$ be  probability measures on $D([0,T];\rrd)$ such that their time marginals,  $\mathbb{P}\circ X^{-1}_t$, $\wt{\mathbb{P}}\circ X^{-1}_t$ have densities  $\call_{X_t}$, and $\wt\call_{X_t}$ respectively, w.r.t. Lebesgue measure for all $t\in[0,T]$ such that
	\begin{equation}\label{e4.4}
		((t,x)\to\call_{X_t}(x)),\ ((t,x)\to\wt\call_{X_t}(x))\in L^\9((0,T)\times\rrd).\end{equation}
If $\mathbb{P}$ and $\wt{\mathbb{P}}$ are solutions to the nonlinear martingale problem $(K_\Box,C^2_c)$, i.e. $($see Definition  {\rm\ref{d4.2})}, they are solutions to the linear martingale problems cor\-res\-pon\-ding to  $(K_{\call_{X_t}},C^2_c)$, $(K_{\call_{\wt X_t}},C^2_c)$, respectively, 
then $\mathbb{P}=\wt{\mathbb{P}}$.
\end{theorem}

\begin{proof} Clearly, by Dynkin's formula, both
$\mu_t(dx):=\call_{X_t}(x)dx$ and $\wt\mu_t(dx):=\call_{\wt X_t}(x)dx,\ t\in[0,T],$ solve the \FP\ equation \eqref{e1.10} with the same initial condition $u_0(dx):=u_0(x)dx$, hence sa\-tisfy \eqref{e3.1} with $y_1(t):=\call_{X_t}$ and $y_2(t):=\call_{\wt X_t}$. Hence, by Theorem \ref{t3.1}, $\call_{X_t}=\call_{\wt X_t},\mbox{ for all }t\ge0,$ since $t\mapsto\call_{X_t}(x)dx$ and $t\mapsto\call_{\wt X_t}(x)dx$ are both narrowly continuous and are probability measures for all $t\ge0$, so both are in $L^\9(0,T;L^1\cap L^\9)\subset L^\9(0,T;L^2)$. 

Now, fixing $\call_{X_t}$ from above, consider the  equation
\begin{eqnarray}
 \label{e4.5}
	v_t+\Psi(-\Delta)\(\dd\frac{\beta(\call_{X_t})}{\call_{X_t}}\,v\)+{\rm div}(Db(u)v)=0,\ \
 	v(0,x)=u_0(x), 
\end{eqnarray}
again in the  distributional sense analogous to \eqref{e1.10}. Then, by Theo\-rem \ref{t3.2} we conclude that $\call_{X_t}$, $t\in[0,T]$, is the unique solution to \eqref{e4.5} in $L^\9(0,T;L^1\cap L^\9)$. Clearly, both $\mathbb{P}$ and $\wt{\mathbb{P}}$ solve the (linear) martingale pro\-blem with   initial condition $u_0(dx):=u_0(x)dx$ corresponding to $(K_{\call_{X_t}},C^2_c)$. Since the above is true for all $u_0\in L^1\cap L^\9$, and also holds when we consider \eqref{e1.1}, resp. \eqref{e4.5}, with start in $s>0$ instead of zero, it follows by exactly the same arguments as in the proof of Lemma 2.12 in \cite{24'} that $\mathbb{P}=\wt{\mathbb{P}}$.\end{proof}

\begin{theorem}\label{t4.4}  For $s\in[0,\9)$ and $\zeta\in\mathcal{Z}:=\{\zeta\equiv\zeta(x)dx\mid\zeta\in L^1\cap L^\9$, $\zeta\ge0,$ $|\zeta|_1=1\}$, let $\mathbb{P}_{(s,\zeta)}$ denote the solution to the nonlinear martingale problem corresponding to $(K_\Box,C^2_c)$ with the initial condition $\zeta$ at the initial time $s$ from Theorems {\rm\ref{t4.1}} and {\rm\ref{t4.3}}. Then, $\mathbb{P}_{(s,\zeta)},$ $(s,\zeta)\in[0,\9)\times\mathcal{Z}$, form a nonlinear Markov process in the sense of Definition {\rm2.1} in {\rm\cite{18a}}, i.e. in the sense of McKean {\rm\cite{17a}}.   
\end{theorem}
\begin{proof} The assertion follows from Corollary 3.8 in \cite{18a} (see also Example (iii) in Section 4.2 of \cite{18a} for the special case with $\Psi(r):=r^s,$ $s\in\(\frac12,1\)$).\end{proof}
\begin{remark}\label{r4.5} \rm Equation \eqref{e4.3} in Theorem \ref{t4.1} says that our solution $u$ of \eqref{e1.1} from Theorem \ref{t2.3} is the one dimensional time marginal law density of a cadl\`ag nonlinear Markov process. This realizes McKean's vision formulated in \cite{17a}   for solutions to nonlinear parabolic  PDEs, namely to identify the solutions to the latter as one-dimensional time marginal law densities of a nonlinear Markov process.
\end{remark}

\section*{Appendix}

\subsection*{A1. Representation and properties of the integral\\\ \hspace*{10mm} kernel   of $(\vp I+\Psi(-\Delta))\1$}\label{appn} 
	
	Let $\vp>0$. We have, for $u\in L^2,\ \xi\in\rrd$,  
\begin{eqnarray*}
	\calf((\vp I+\Psi(-\Delta))\1u)(\xi)\!\!\!\!\! 
	&=&\!\!\! \frac1{\vp+\Psi(|\xi|^2)}\ \calf u(\xi) 
	 = \int^\9_0 e^{-\vp t}\ e^{-t\Psi(|\xi|^2)}dt\,\calf u(\xi)\\
	&\buildrel{(1.7)}\over{=\!=\!}&\!\!\!\int^\9_0\ e^{-\vp t}\int^\9_0 e^{-r|\xi|^2}\eta^\Psi_t(dr)dt\,\calf u(\xi)\\
	&=&\!\!\!(2\pi)^{\frac d2}\dd\int^\9_0 e^{-\vp t}\int^\9_0\calf(p_r)(\xi)\eta^\Psi_t(dr)dt\,\calf u(\xi),\end{eqnarray*} 
	$$p_r(x):=\frac1{(4\pi r)^{\frac d2}}\ e^{-\frac1{4r}\,|x|^2},\ x\in\rrd.$$ 
\n	Hence, defining
	$$g^\Psi_\vp(x):=
	\int^\9_0 e^{-\vp t}\int^\9_0 p_r(x)\eta^\Psi_t(dr)dt,\ x\in\rrd,\eqno{\rm(A.1)}$$we have
	$$(\vp I+\Psi(-\Delta))\1 u=g^\Psi_\vp*u.$$
	Since $\eta^\Psi_t,\ t\ge0,$ are probability measures, we have
	$$\vp\int_\rrd g^\Psi_\vp  dx=1,\ \ \ff\vp>0.\eqno{\rm(A.2)}$$

\subsection*{A2. The uniqueness of equation \eqref{2.27'''}}

\mk\n{\bf Lemma A.} {\it Assume that $d\ge2$ and let $y_1(\lbb),\ y_2(\lbb)\in L^1\cap L^\9$ be two distributional solutions to \eqref{2.27'''}. Then, if Hypotheses {\rm(i)--(iii)} hold, there exist $\wt\lbb_0\in(0,\lbb_0)$ such that for all $\lbb\in(0,\wt\lbb_0)$ we have $y_1(\lbb)=y_1(\lbb).$}

\begin{proof} The proof of this lemma is completely analogous to that of Lemma A in \cite{a9}. Therefore, we omit it here.
\end{proof}

\n{\bf Acknowledgments.} This work was supported by the Deutsche Forschungsgemeinschaft (DFG, German Research Foundation) -- Project ID 317210226 -- SFB 1283 and by a grant of the Ministry of Research, Innovation and Digitization, CNCS--UEFISCDI project PN-III-P4-PCE-2021-0006, within PNCDI~III. 
A part of this work was done during very pleasant stays of the third named author at the University of Madeira as a guest of the second named author.


\begin{thebibliography}{nn}
\bibitem{1} Barbu, V., \textit{Nonlinear Differential Equations of Monotone Type in Banach Spaces},  Berlin. Heidelberg. New York, Springer, 2010.\vspace*{-2,5mm}
	
\bibitem{a2''} Barbu, V., R\"ockner, M.,   From \FP\ equations to solutions of distribution dependent SDE, \textit{Ann. of Probab.},   {48} (2020),   1902--1920.\vspace*{-2,5mm}
	
\bibitem{a2'''} Barbu, V., R\"ockner, M.,   Solutions for nonlinear \FP\ equations with measures as initial data and McKean-Vlasov equations.  \textit{J.~Funct. Anal.}, {280} (7) (2021),  1--35.\vspace*{-2,5mm} 
	
\bibitem{2} Barbu, V., R\"ockner, M.,   The evolution to equilibrium of solutions to nonlinear Fokker-Planck equations, \textit{Indiana Univ. Math.~J.},  {72} (1) (2023),  89--131.\vspace*{-2,5mm} 	
	
\bibitem{4} Barbu, V., R\"ockner, M.,  Uniqueness for nonlinear  Fokker-Planck equations and for McKean--Vlasov SDEs: the degenerate case, \textit{J. Funct. Anal.},  {285} (4) (2023).\vspace*{-2,5mm} 
	
\bibitem{a9} Barbu, V., R\"ockner, M., Nonlinear \FP\ equations with fractional Laplacian and McKean--Vlasov SDEs with L\'evy noise, {\it Probab. Theory Relat. Fields}, 189 (2024) 849-878.\vspace*{-2,5mm}   

\bibitem{6b} Barbu, V., R\"ockner, M., {\it Nonlinear \FP\ Flows and their Probabilistic Counterparts}, Lecture Notes in Mathematics, vol. 2353, Springer, 2024.\vspace*{-2,5mm}
	
\bibitem{5} Brezis, H., Crandall, M.G.,  Uniqueness of solutions of the initial-value problem for $u_t-\Delta\beta(u)=0$,  \textit{J. Math. Pures et Appl.},   {58} (1979),  153--163.\vspace*{-2,5mm} 
	
\bibitem{8a} Carrillo, J.A.,  Entropy solutions for nonlinear degenerate pro\-blems, \textit{Archives Rat. Mech. Anal.}, {147} (1999),   269--361.\vspace*{-2,5mm} 
	
\bibitem{a9'} Carmona, R., Delarue, F.,  \textit{Probabilistic Theory of Mean Field Games with Applications}, I--II. Springer, 2017.\vspace*{-2,5mm} 
	
\bibitem{9prim} del Teso, F., Endal, J., Jakobsen, E.R., Uniqueness and pro\-per\-ties of distributional solutions of nonlocal equations of porous media type, \textit{Adv. Math.}, {305} (2017), 78-143.\vspace*{-2,5mm}
	
\bibitem{7} De Pablo, A., Quiros, F., Rodriguez, A., Vasquez, J.L.,  A general fractional porous medium equation, \textit{Comm. Pure Appl. Math.}, {65} (9) (2012),  1242--1284.\vspace*{-2,5mm} 
	
\bibitem{7a} De Pablo, A., Quiros, F., Rodriguez, A., Vasquez, J.L.,   A fractional porous medium equation. \textit{Adv. Math.}, {226} (2) (2010),   1378--1409.\vspace*{-2,5mm} 
	
\bibitem{b17'} Figalli, A.,  Existence and uniqueness of martingale solutions for SDEs with rough or degenerate coefficients. \textit{J. Funct. Anal.}, {254} (1) (2008),   109--153.\vspace*{-2,5mm} 

\bibitem{16'} Fukushima, M., Oshima, Y., Takeda, M.,   \textit{Dirichlet Forms and Sym\-me\-tric Markov Processes}, de Gruyter, 2011,  x+489 pp.\vspace*{-2,5mm} 
	
\bibitem{a11'} Funaki, T.,  A certain class of diffusion processes associated with nonlinear parabolic equations, \textit{Z. Wahrsch. Verw. Gebiete},  {67} (3) (1984),   331--348.\vspace*{-2,5mm} 
	
\bibitem{MR92} Ma, Zhi Ming, R\"ockner, M.,  Introduction to the theory of (non\-symmetric)  Dirichlet forms, \textit{Universitext}. Berlin, Springer Verlag, 1992, vi+209~pp.\vspace*{-2,5mm} 
	
\bibitem{17a} McKean, H.P.,  A class of Markov processes associated with nonlinear 	parabolic equations, \textit{Proc. Nat. Acad. Sci. U.S.A.}, {56}  (1966),   1907--1911.\vspace*{-2,5mm} 
	
\bibitem{a15'} Pierre, M.,  Uniqueness of the solutions of $u_t-\Delta\varphi(u)=0$ with initial data measure, \textit{Nonlinear Anal.}, {6} (2) (1982),    175--187.\vspace*{-2,5mm} 
	
\bibitem{18a}  Rehmeier, M., R\"ockner, On nonlinear Markov processes in the sense of McKean,  {\it Journal of Theoretical Probab.}, vol. 38 (2025), 36 pp.\vspace*{-2,5mm}
	
\bibitem{11''} Ren, P., R\"ockner, M., Wang, F.Y.,  Linearization of nonlinear \FP\ equations and applications, \textit{J. Differential Equations}, {322} (2022),   1--37.\vspace*{-2,5mm} 
	
\bibitem{23'} R\"ockner, M., Wu, W.,  Xie, Y.,   Stochastic generalized porous media equations over $\sigma$-finite measure spaces with non-continuous diffusivity function, {\it Potential Analysis}, 61 (2024), 731--773.\vspace*{-2,5mm} 
    	
\bibitem{7aa} R\"ockner, M., Xie, L., Zhang, X., Superposition principle for nonlocal \FP--Kolmogorov operators,   \textit{Probab. Theory Rel. Fields},   {178} (3--4) (2020),  699--733.\vspace*{-2,5mm} 
	
\bibitem{12''} Schilling, R., Song, R., Vondra\v cek, Z.,   \textit{Bernstein Functions}, de Gruyter, 2012.\vspace*{-2,5mm}  	
	
\bibitem{24'} Trevisan, D.,  Well-posedness of multidimensional diffusion processes with weakly differentiable coefficients, \textit{Electron J. Probab.}, {21}  (2016),   22--41.\vspace*{-2,5mm} 
	
\bibitem{11'} 	V\'asquez, J.L.,  Nonlinear diffusion with fractional Laplacian ope\-ra\-tors. In \textit{Nonlinear partial differential equations. The Abel Symposium 2010}. pp. 271--298. Abel Symposia, {7}. Berlin. Heidelberg, (2012).   doi:10.10007/978-3-642-25361-4\_15.\vspace*{-2,5mm} 
	
\end{thebibliography}
\end{document}